\documentclass[]{interact}

\usepackage{epstopdf}
\usepackage{hyperref}
\hypersetup{
    colorlinks=true,
    linkcolor=blue,
    citecolor = blue,
    filecolor=black,      
    urlcolor=blue,
}
\usepackage[caption=false]{subfig}

\usepackage[numbers,sort&compress]{natbib}
\bibpunct[, ]{[}{]}{,}{n}{,}{,}
\renewcommand\bibfont{\fontsize{10}{12}\selectfont}

\theoremstyle{plain}

\theoremstyle{definition}

\theoremstyle{remark}

\begin{document}

\articletype{}

\title{Structure from Appearance: Topology with Shapes, without Points}

\author{
\name{Alexandros Haridis\textsuperscript{a}\thanks{CONTACT A. Haridis. Email: charidis@mit.edu} 
}
\affil{\textsuperscript{a}Department of Architecture, Massachusetts Institute of Technology, 77 Massachusetts Ave, Room 10-303, Cambridge, MA 02116, USA
}
}

\maketitle

\begin{abstract}
A new methodological approach for the study of topology for shapes made of arrangements of lines, planes or solids is presented. Topologies for shapes are traditionally built on the classical theory of point-sets. In this paper, topologies are built with shapes which are formalized without points and with structures defined from their parts. An interpretative, aesthetic dimension is introduced according to which the topological structure of a shape is not inherited from an ambient space but is induced based on how it is interpreted into sets of parts. The proposed approach provides a more natural, spatial framework for studies on the mathematical structure of design objects, in art and design. More generally, it shows how mathematical constructs (here, topology) can be built directly in terms of objects of art and design, as opposed to the more common opposite approach where objects of art and design are subjugated to canonical mathematical constructs.

\end{abstract}


\section{Introduction}

This paper is about a particular type of structure, namely about \emph{topological structure} or \emph{topology}, and it is concerned with how this structure can be defined and studied in terms of interpretations of the appearance of \emph{shapes}. Appearance and interpretation are two concepts that are in many ways foreign to the development of topology in mathematics. Both are, however, inherent in art and design, in the ways designers and artists work when they explore ideas about form and composition. Interpretation is understood here in an aesthetic sense: it is the intuitive and natural act we commonly find in art and design, of describing the sensuous ``surface" of a design object or an artwork---its appearance (``how it looks")---according to parts we choose to perceive or make salient in it\footnote{``Part" and ``interpretation" are concepts that appear in various formal systems; e.g. \cite{LeonardGoodman40, Simons87}. In this paper, they are understood pictorially, and are formalized in a particular formal system for shapes (\emph{Section 2}).}. The shapes I have in mind are the pictorial ones, formed in arrangements of lines or planes, but also the shapes of physical things, formed in arrangements of solids. Such shapes are the basis for drawings, sketches, models, visual compositions and other means of creative expression. My goal in this paper is to study topology directly with such shapes and to make the intuitive act of interpreting their appearances into parts an actual ``mechanism" for inducing or generating topologies.

Traditionally, the objects upon which one wishes to confer a topology must be represented as sets, in one space or another; the generality of topological constructions rests on this fact \cite{Kurat72}. In general topology (point-set topology), topologies for shapes are modeled after the classical concept of ``topological space" \cite{Munkr00}: a shape is considered a subspace of an underlying ambient space (e.g., Euclidean space), which possesses a set of points \emph{and} a predefined system of subsets containing the points; the subsets determine ``relations between points" (Figure 1). The point-set view of shapes is convenient for most applications. It allows one to transfer constructions from the existing literature of topology in mathematics to the study of topology for shapes, without significant alterations. This is the approach followed in many areas of engineering design \cite{RosenPeters96}. Most notably, it is the basis for geometric modeling in computer graphics and Computer-Aided Design \cite{Requi77}.

\begin{figure}[ht]
\centering
\includegraphics{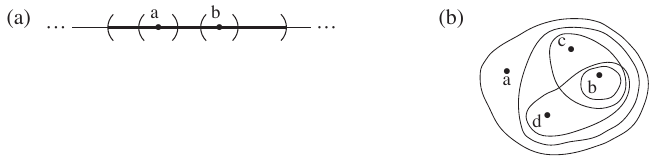}
\caption{Graphical representation of (a) infinite and (b) finite topological spaces. A topological space possesses a set of points and a system of ``open sets" containing the points.}
\label{Figure1}
\end{figure}

The point-set view of shapes, however, has been characterized as unintuitive and unnatural for domains where there is a strong aesthetic component, such as architecture, design, or the visual arts (evidence from empirical studies can be found in \cite{VerstijnenHennessey98, StonesCassidy2010}; a more general, formal argument is established in \cite{JowEarlStiny2019}). In art and design, shapes are understood and manipulated as geometric objects in their own right---they are not understood in terms of other things (e.g., in terms of specialized point-sets). For design purposes, shapes are synthesized from elements that have an actual appearance and ``physical" presentation, with no apparent subdivision scheme (sets of ``infinitesimally small" elements of trivial or no content/extension are not naturally compatible with this). Drawings and models, constructed out of shapes, must be seen or touched to be appreciated and evaluated, as wholes or in parts, by designers and artists.

Point-set topology, even though broadly applied to various design and engineering contexts, is not necessarily meant to address spatial or aesthetic concerns. Some aspects of this incompatibility, particularly related to boundaries of shapes, are examined in \cite{Earl97}. With respect to appearance, the idea that a shape can obtain structures from different ways of interpreting its appearance into parts (the designer's or artist's ``way of working") is essentially absent. Topological investigations concerning shapes usually go in the opposite direction. They establish similarities, or correspondences, between the topologies that may underlie seemingly different shapes---so that, topologically similar but perceptually different appearances can be reduced into the same ``structural class". Classificatory studies concerning knots, braids and links \cite{Rolfsen2003}, surfaces and simplicial complexes \cite{Lefsch42, Munkr00}, are a few examples of this general direction in the mathematical literature (we can detect the pursuit for topological similarity even in popular opinion, in phrases like ``a torus is topologically the surface of a donut, or a rubber tire, or a coffee mug"). In \cite{Vard2017}, the reader can find a historical examination (and confirmation) of the absence of appearance in structural mathematics in general and how it has influenced morphological studies in architecture and design.

With this paper I aim to show that pictorial/spatial ideas and ``ways of working" coming from non-mathematical areas with a strong aesthetic component (architecture, design or the visual arts), can actually acquire mathematical substance. This suggests a way of doing mathematics within art and design according to which the relevant mathematical constructs (here, topology) are built in terms of objects directly inspired and taken from art and design. As opposed to working in the more common opposite direction, where those objects have to be first subjugated to a canonical mathematical scheme or format.

The basis of this investigation is a formalization of shapes that comes from the theory of computation in design developed around the shape grammar formalism \cite{Stiny75, Stiny2006}. Shape grammars were originally invented as a model of computation for describing languages of shapes with certain stylistic interest (geometric paintings, sculptures, patterns, architectural plans, ornamentation \cite{StinyGips72, Knight94}). They are now the basis for a full-fledged theory of calculating in art and design.

The mathematical theory of shapes behind shape grammars provides an alternative to the point-set view of shapes. It is captured in the \emph{algebras of shapes} $U_i$, with a standard classification of shapes in terms of the dimensionality $i$ of their basic elements \cite{Stiny75, Stiny91}\footnote{In the literature, the notation $U_{ij}$ is sometimes used to specify the dimensionality $j$ of the ``space" in which a shape is formed. For example, arrangements of lines in two dimensions are shapes in an algebra $U_{12}$. The index $j$ does not play a role in the context of this paper and is thus omitted (it is assumed of course that $i \leq j$).}. Technical details are both intuitive from a design and artistic point of view and formally precise in the mathematical sense. In principle, the algebras $U_i$ provide a point-free theory of shape (when lines, planes or solids are involved), equipped with a primitive concept of \emph{part} \emph{relation}. The part relation puts together ``appearance" and ``interpretation" at the center of the theory. We essentially have a mechanism that fixes mathematically an intuitive notion of interpretation: that shapes have parts, which are shapes embedded in them, and that those parts are not fixed and given (as members of a set or points of a topological space are given), but are ambiguous and thus can be seen or distinguished in many ways. 

In the literature of shape grammars, the implications of using the algebras $U_i$ as a basis for calculation have been examined in depth (the most recent summary of this work is in \cite{Stiny2006}). Research related to the structure of shapes in algebras $U_i$ has also appeared although in fewer numbers \cite{Knight88, Stiny94, Krst96, JowEarlStiny2019}. Most of this research concentrates on spotlighting how structure intervenes when designers use Computer-Aided Design systems, and in what ways it discourages, or eliminates, fluid, open-ended manipulations of a shape's parts (e.g., \cite{Stiny94, JowEarlStiny2019}). In some publications, structures for shapes are examined in more mathematical detail \cite{Krst96}. Specifically for the subject of topology, however, previous work (mainly in \cite{Earl97}, \cite{Krst96} and \cite{Stiny94}) has not examined in full detail the possibility of using the algebras $U_i$ as an alternative foundation to the classical point-set view of shapes for doing topology. My purpose here is to contribute to this end.

In this paper, I develop a framework for topology with shapes that takes the algebras $U_i$ as the underlying theory of shape. In this framework, an interpretative, aesthetic dimension is introduced in the sense that a topological structure is not inherited from a predefined ambient space, but it is \emph{induced} based on how a shape is interpreted into finite sets of parts. \emph{Basis}, \emph{continuity}, \emph{connectedness}, and other foundational topological concepts, are formulated directly in terms of shape and topological structures defined by parts. This paper focuses on \emph{finite topological structures}, that is to say, topologies with finitely many parts---these topologies can be understood as finitistic descriptions of shapes and are natural tools for describing designs. In the absence of a point-set substrate (a desirable effect of the algebras $U_i$), topological concepts normally defined in terms of points are now formulated only in terms of part relations, and their formulation is driven by spatial intuitions. In \emph{Section 2}, I present the algebras of shapes $U_i$, the background required for this work. The technical material on finite topology for shapes is in \emph{Sections 3} through \emph{7}.

I do not specifically focus on how theorems and results for topological spaces can be mimicked or generalized in this framework. In certain cases, it is a straightforward task; some evidence and discussion is provided. More interestingly though, in light of a point-free framework for topology and the introduction of appearance as a key factor for determining structure, I show how concepts that have been worked out for topological spaces do not always transfer naturally. Sets (or spaces) with definite subdivision into points are not categorically the same objects as shapes without definite parts. This change in the underlying object of study leads to modifications and adaptations that many times require some ingenuity to be meaningful.

The technical material here can also serve as a first introduction to point-free thinking in topology if one wishes to do so through shapes (see \cite{Johnstone2001}, \cite{Menger40}, and \cite{PicPultr2010}, for general references to point-free topology in mathematics). At least when it comes to the finitistic case, the paper provides a much more intuitive and actually \emph{spatial} framework, compared to other point-free frameworks available in the mathematical literature (such as the theory of locales and frames, e.g., \cite{PicPultr2010}). In the final section of this paper, I comment on how this paper contributes to studies of the mathematical structure of design objects, in art and design.

\section{Background: shapes and visual calculating}

\subsection{Shapes and algebras $U_i$}

Traditionally, the beginning point of topology is the classical theory of point-sets (e.g., \cite{Lefsch42, Munkr00, Barm2011}). In this paper, topology is built on the theory of shape that originates in the shape grammar formalism \cite{Stiny75, Stiny2006}. In this section, I treat those topics of the theory, and only those, that will be needed in \emph{Sections 3} through \emph{7} of this paper. I provide an elementary and rather informal overview, with relatively few technical arguments. For further background and more technical coverages, see \cite{Krst96, Stiny75, Stiny92, Stouffs94}.

In the shape grammar formalism, shapes are understood in a manner analogous to how artists, designers or architects create, appreciate and manipulate drawings and physical models in practice, especially in the early stages of the creative process. Thus, the theory of shape is driven by intuitions coming directly from art and design.

\begin{figure}[ht]
\centering
\includegraphics{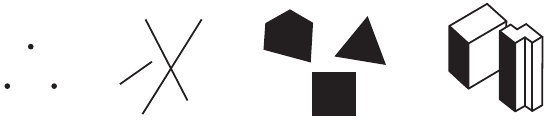}
\caption{Basic elements for shapes in algebras $U_i$: points, lines, planes and solids.}
\label{Figure2}
\end{figure}

Shapes have finite spatial extension and are made up of basic elements: points, lines, planes and solids. These are things like those in Figure 2. The first three---points, lines, and planes---are pictorial elements, and we can easily draw them on a sheet of paper. Solids represent the shapes of real, physical things. We can visualize solids in a drawing, as in the drawing in Figure 2, but the drawing itself is made with lines and planes (not solids). Shapes made with lines, planes or solids are unanalyzed and their parts can be seen in many ways. Though shapes, in general, can be made from any combination of basic elements, either of a single kind or of different kinds, the investigations in this paper are concerned only with shapes made with basic elements of a single kind.

Shapes are formalized in algebras $U_i$, originally invented for purposes of calculation with shape grammars \cite{Stiny91}. In the algebras, shapes are classified in terms of the dimensionality $i$ of their basic elements: $i$ = 0 for points, $i$ = 1 for lines, $i$ = 2 for planes, and $i$ = 3 for solids. Some properties of those basic elements are summarized in Table 1. Other basic elements, for example, curves or surfaces, can also be considered to extend the algebras $U_i$ (\cite{Stiny2006, JowEarl2015} provide further elaboration).

\begin{table}[ht]
\tbl{Properties of basic elements.}
{\begin{tabular}{lcccc} \toprule \\
 Basic element & Dimension & Boundary & Content & Part relation \\ \midrule
 Point & 0 & none & none & identity \\
 Line & 1 & two points & length & partial order \\
 Plane & 2 & three or more lines & area & partial order \\ 
 Solid & 3 & four or more planes & volume & partial order \\\bottomrule
\end{tabular}}
\label{basicelements-table}
\end{table}

Three things define an algebra of shapes. First, we have the shapes themselves. A shape is a finite arrangement (a set) of basic elements of a certain kind, which are \emph{maximal} with respect to each other. Second, an algebra is equipped with a \emph{part relation} denoted with $\leq$ that drives the operations of sum (+), product ($\cdot$) and difference (-) for shapes and enables recognition of parts in any given shape. And third, we have \emph{transformations} for changing a given shape into others. A standard example are the Euclidean transformations, but the algebras allow linear transformations and possibly other kinds, too. I focus on the first two of these characteristics of the algebras, for they underlie all topological constructions developed later on. 

A shape is uniquely defined by the set of its maximal elements \cite{Stiny75}:

\vspace{0.1in}
\noindent \emph{The smallest set containing the biggest basic elements that combine to form the shape}.

\vspace{0.1in}
\noindent Given a shape, one can follow a precise algorithmic procedure to reduce the shape to its maximal elements. The interested reader can refer to \cite{Stiny75, Stiny2006} for technical details. An intuitive way to understand maximality is by analogy to how skilled drafters make line-drawings, either manually or with the help of a digital computer. To make a drawing, a drafter draws the ``biggest lines" containing all lines of the drawing (this is a standard practice for students in design schools). The \emph{fewest} such lines needed to complete the drawing correspond precisely to the maximal elements of the drawing. Figure 3 shows examples of shapes made with points, lines and planes, and their corresponding maximal elements.

\begin{figure}[ht]
\centering
\includegraphics[scale=0.92]{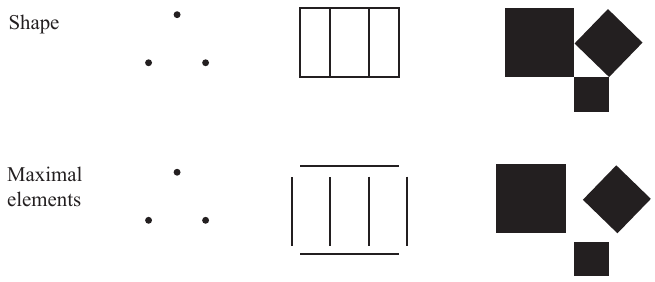}
\caption{Shapes made with points, lines and planes, and their corresponding maximal elements.}
\label{Figure3}
\end{figure}

The maximal representation of a shape in the algebras is for purposes of finite describability. This augments the computational potential. The parts of a shape, however, are not only those that are explicitly defined by maximal elements but also those that can simply be seen or recognized in them figuratively. The part relation ($\leq$) is the main relation that enables this, for all basic elements in the algebras $U_i$.

The only point that is part of another point is the point itself. A point is like a drawing that has no proper nonempty parts, regardless of how the point looks (see also \cite{Haridis2020}). For shapes made with points in algebra $U_0$, the part relation behaves as an \emph{identity relation} (an \emph{equivalence}). One can explicitly enumerate the parts of a shape by forming sets of zero, one or more points; this is illustrated in Figure 4a, for a shape made with three points.

\begin{figure}[ht]
\centering
\includegraphics[width=\textwidth]{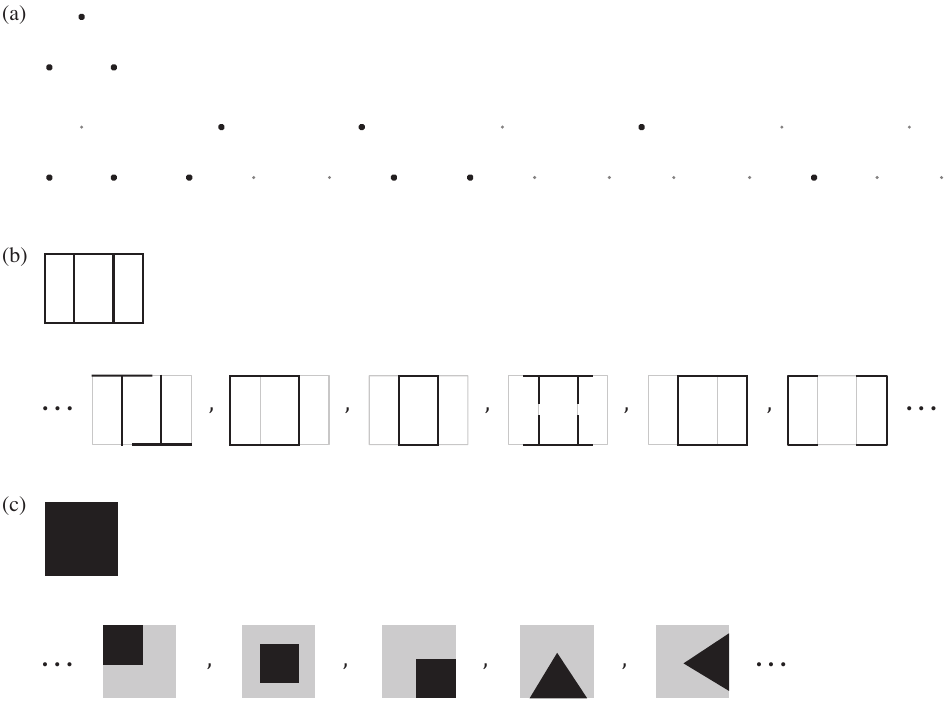}
\caption{(a) A shape made with three points and enumeration of its parts. (b), (c) Shapes made with lines or planes have uncountably many parts in them.}
\label{Figure4}
\end{figure}

In an algebra $U_0$, a shape and its parts form a finite Boolean algebra; the atoms are the individual points of the shape. For example, in Figure 4a, every part comes with its complement relative to the shape, and the overall set of parts is closed under finite sums and products (the operations of sum, product and relative complement for shapes are explained a little later). This algebra has for the bottom element a special shape called the \emph{empty shape}, denoted with \emph{0}, and for the top element the shape itself. The \emph{empty shape} (at times called \emph{empty part}, depending on the context), is a shape with no parts; intuitively, it corresponds to a blank sheet of paper, or an empty drawing. The Boolean algebra formed by a shape and its parts in $U_0$, corresponds to the Boolean algebra formed by a finite set and its subsets.

Shapes in an algebra $U_i$ when $i > 0$ are made up of basic elements that have uncountably many basic elements of the same kind embedded in them. That is to say, lines have parts that are lines of finite but nonzero length (e.g., Figure 5), planes have parts that are planes of finite but nonzero area, etc. The content (measure) of a basic element, other than a point, is never zero (nor is the content of its nonempty parts). For shapes in an algebra $U_i$ where $i > 0$, the part relation behaves as a \emph{partial order} (see \cite{Stiny75} for the mathematical argument). For any given shape one can \emph{trace}, or equivalently \emph{recognize} or \emph{see}, uncountably many different parts embedded in it (including, but not limited to, its maximal parts), which are shapes made up of basic elements of the same kind as the shape itself. A graphical illustration is given in Figure 4b and 4c.

\begin{figure}[ht]
\centering
\includegraphics{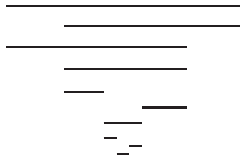}
\caption{The parts of a line in algebra $U_1$ are never points but lines of finite but nonzero length.}
\label{Figure5}
\end{figure}

In an algebra $U_i$ where $i > 0$, a shape and its parts form an infinite Boolean algebra. Every part of the shape has a complement relative to the shape, finite sums and products of parts determine parts of the same shape, the \emph{empty shape} is the bottom element and the shape itself is the top element. This algebra, however, is atomless. It is also not complete, because infinite sums or products of parts do not determine parts/shapes necessarily \cite[p.~208]{Stiny2006}.

Even though shapes have uncountably many parts to see in them (when $i > 0$), they are not themselves infinite objects: a shape does not come equipped with the set of its parts given all at once! Instead, a shape is a \emph{finite} object, kept in maximal element representation, but has a natural ability to be divisible into parts indefinitely (real physical drawings and models in art and design have these properties, too). Intuitively, you choose parts by tracing them or equivalently seeing them. There may be uncountably many to see, but each time you see finitely many of them---you do not see infinitely many parts all at once. This distinguishing characteristic of shapes, namely their finite describability and at the same time their potential for indefinite divisibility, is central to the topological constructions developed in this paper. A few other background concepts are needed before we proceed.

All basic elements besides points have a unique boundary (Table 1). The boundary of a basic element is made up of basic elements of exactly one dimension lower than itself (e.g. the boundary of a line is formed by two points). This leads to the following result: the boundary of a basic element is not a part of the basic element. More generally, given a shape in an algebra $U_i$, for $i > 0$: (i) the boundary of the shape is a \emph{shape} made with basic elements of dimension $i$ - 1, and (ii) this boundary shape is not a part of the original shape. (For further details on shapes and their boundaries, and a comparison with the corresponding concepts in point-set topology, see \cite{Earl97}.)

The algebras $U_i$ are closed under operations of sum (+), product ($\cdot$), and difference (-). These operations mimic very closely how drawings are created and manipulated in practice, in art and design; their technical details are covered in \cite{Stiny75, Krishn92a, Krishn92b, Stiny2006}.

Suppose $S$ and $S'$ are two shapes, both made with the same kind of basic elements (e.g., $S$ and $S'$ are both made with points or both made with lines). Say that $S$ is part of $S'$, and denote it with $S \leq S'$, if every maximal element of the first is embedded in a maximal element of the second (more generally, under an appropriate transformation). An intuitive way to understand this is by an analogy to manual drawing: first draw the shape $S$ and then the shape $S'$; if the resulting drawing is $S'$, then $S$ is part of $S'$. 

A brief description of sum, product and difference for shapes now follows.

\vspace{12pt}

{\small \noindent \emph{Sum}\;\; $S$ + $S'$, is the unique shape formed by adding two shapes together; it corresponds to the act of drawing shape $S$ and then shape $S'$. The resulting shape satisfies two conditions: (i) $S$ and $S'$ are part of the sum, and every part of the sum has a part that is part of one shape or the other. If $S \leq S'$, then $S$ + $S'$ = $S'$. We also have that $S$ + $S$ = $S$, for any shape $S$.}

\vspace{12pt}

{\small \noindent \emph{Difference}\;\; $S$ - $S'$, is the unique shape formed by subtracting from $S$ all parts that are shared with $S'$; it corresponds to the manual act of erasing parts from $S$. Every part of the difference is a part of $S$ but not $S'$.}

\vspace{12pt}

{\small \noindent \emph{Product}\;\; $S \cdot S'$, is the largest part shared with both $S$ and $S'$. Alternatively, the product is formed by the difference: $S \cdot S'$ = $S$ - ($S$ - $S'$). When two shapes share no parts, $S \cdot S'$ is equal to the \emph{empty shape}. In this case, it also follows that $S$ - $S'$ = $S$ and $S'$ - $S$ = $S'$. We also have that $S \cdot S$ = $S$, for any shape $S$.}

\vspace{12pt}

\begin{table}[ht]
\tbl{Operations of sum, product and difference for shapes.}
{\begin{tabular}{ccccc} \toprule \\
 $S$ & $S'$ & Sum & Product & Difference \\ \midrule\\
 \includegraphics[scale=0.85]{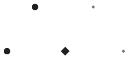} & \includegraphics[scale=0.85]{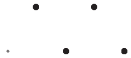} & \includegraphics[scale=0.85]{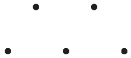} & \includegraphics[scale=0.85]{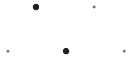} & \includegraphics[scale=0.85]{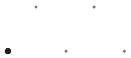} \\\\\\
 \includegraphics[scale=0.95]{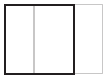} & \includegraphics[scale=0.95]{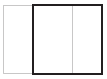} & \includegraphics[scale=0.95]{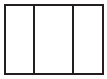} & \includegraphics[scale=0.95]{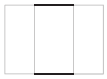} & \includegraphics[scale=0.95]{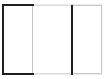} \\\\\\
 \includegraphics[scale=0.95]{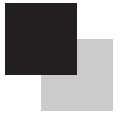} & \includegraphics[scale=0.95]{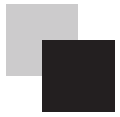} & \includegraphics[scale=0.95]{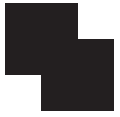} & \includegraphics[scale=0.95]{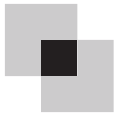} & \includegraphics[scale=0.95]{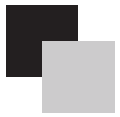} 
 \\\\\bottomrule
\end{tabular}}
\label{operations-table}
\end{table}

An illustration of the three operations is given in Table 2, for shapes made with points, lines and planes. Notice that for shapes made with points, sum, product and difference coincide, respectively, with union, intersection and difference for sets.

Last, shapes have no absolute complements because the algebras $U_i$ do not have a unit \cite[pp.~205--209]{Stiny2006}, i.e., there is no largest shape that contains all other (such an all-encompassing concept would be meaningless for pictorial and spatial objects in design). A shape can have a complement only \emph{relative to} a specific shape. The relative complement of shape $S$ with respect to shape $S'$ is equal to the shape $S'$ - $S$.

\subsection{Aesthetic interpretation of shapes (and artworks)}

The part relation in the algebras $U_i$ is the driving force of calculations with shape grammars. More than this, the part relation is the formal, mathematical mechanism that enables the possibility of interpreting the appearance of a shape aesthetically.

Previously, it was mentioned that any shape is defined in terms of its maximal elements. Maximal elements, however, are not the only parts that can be recognized or perceived in a shape (whether individually or in combinations). The part relation for basic elements of dimension $i > 0$ supports seeing beyond them. The following is a consequence of that:

\vspace{0.1in}
\centerline{\emph{How a shape is defined is independent of how it is interpreted aesthetically}.}

\vspace{0.1in}
\noindent Consider the shape made with lines in Figure 6a. The maximal elements of this shape are shown in Figure 6b; there are eight lines. Independently of these maximal elements, one can see within the same shape plenty of others. For example, we can see the shape in Figure 6c (this is drawn after Wassily Kandinsky \cite[p.~138]{Kandinsky47}): each line of the shape on the right is part of a maximal element of the shape on the left. We can also see all the shapes in Figure 7; these can be embedded (seen) in the same shape in multiple ways, for example, as in Figure 8. Since basic elements of dimension greater than zero can be divided in any number of ways into any number of parts, there are indefinitely many \emph{different} parts to recognize in any given shape---a graphical example of this is in Figure 9.

\begin{figure}[ht]
\centering
\includegraphics{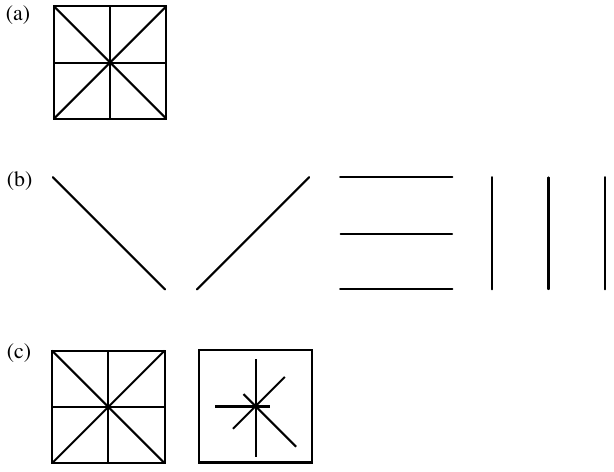}
\caption{(a) A shape made with lines and (b) the maximal elements of this shape. (c) An interpretation of the appearance of the shape in (a) that results in a new shape that is part of the original one.}
\label{Figure6}
\end{figure}

These examples are suggestive of a more general approach to the interpretation of the appearance of a shape that is based on seeing. Interpretation, as it is understood here, is the act of describing the sensuous ``surface" of a shape: whatever it is that it appears to be doing on its surface. It is about \emph{aesthetic interpretation}, that is to say, an interpretation of the way the shape looks. An aesthetic interpretation can be understood as a description of the appearance of a shape in terms of certain parts---the parts that you choose to see in it. For any shape without points, there can be indefinitely many \emph{different} aesthetic interpretations of its appearance, regardless of how this shape was put together to begin with. To say this differently: 

\vspace{0.1in}
\centerline{\emph{The parts used to make the shape are not necessarily the parts you choose to see in it}.}
\vspace{0.1in}

\begin{figure}[ht]
\centering
\includegraphics{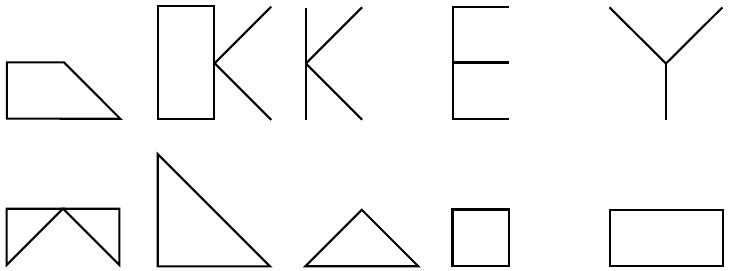}
\caption{Examples of parts embedded in the shape in Figure 6a.}
\label{Figure7}
\end{figure}

\begin{figure}[ht]
\centering
\includegraphics{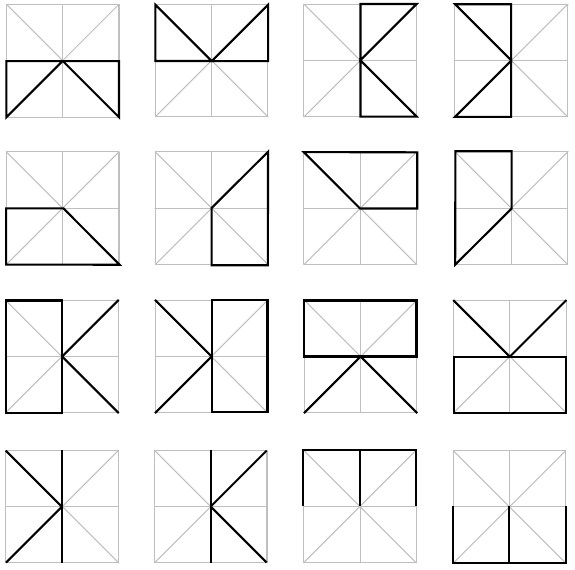}
\caption{Alternative ways of embedding shapes from Figure 7 in the same shape in Figure 6a. These embeddings happen under different Euclidean transformations.}
\label{Figure8}
\end{figure}

\noindent This naturally connects with how artworks, especially the pictorial ones, are appreciated by their different viewers. Things like drawings, sketches, paintings, visual compositions etc., are pictorial objects made with shapes the appearance of which can be appreciated in terms of the parts that a viewer perceives in them. For example, the drawing in Figure 10 by Paul Klee, is an arrangement of lines and curves. Like shapes, the appearance of this drawing can be interpreted in multiple ways, perhaps in terms of the series of parts in Figure 11. Similarly, the drawing in Figure 12 by Pablo Picasso can be considered as a composition of lines or a composition of planes, with a possible interpretation in terms of the parts shown in Figure 13.

\begin{figure}[!ht]
\centering
\includegraphics{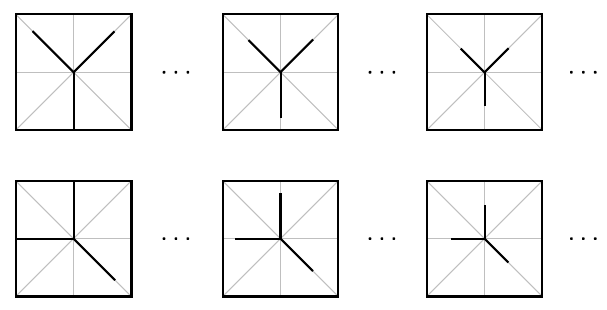}
\caption{Infinitely many parts of the shape in Figure 6a obtained by successive scalings of a part that looks like a capital letter Y (\textit{Source}: Redrawn from Stiny (2006)).}
\label{Figure9}
\end{figure}

These two examples are indicative of what is possible. In both of them, the visible surface of the artwork has been treated as a shape without points. Interestingly enough, interpretations of this kind can happen in an open ended manner for pretty much any aesthetic object (drawing, painting, sculpture, building, etc.), independently of how this object was conceived or made by its original artist. The parts that a viewer appreciates each time in an aesthetic object---the same viewer multiple times, or multiple viewers at the same time---are the parts that the viewer chooses to embed or recognize in it.

\begin{figure}[ht]
\centering
\includegraphics{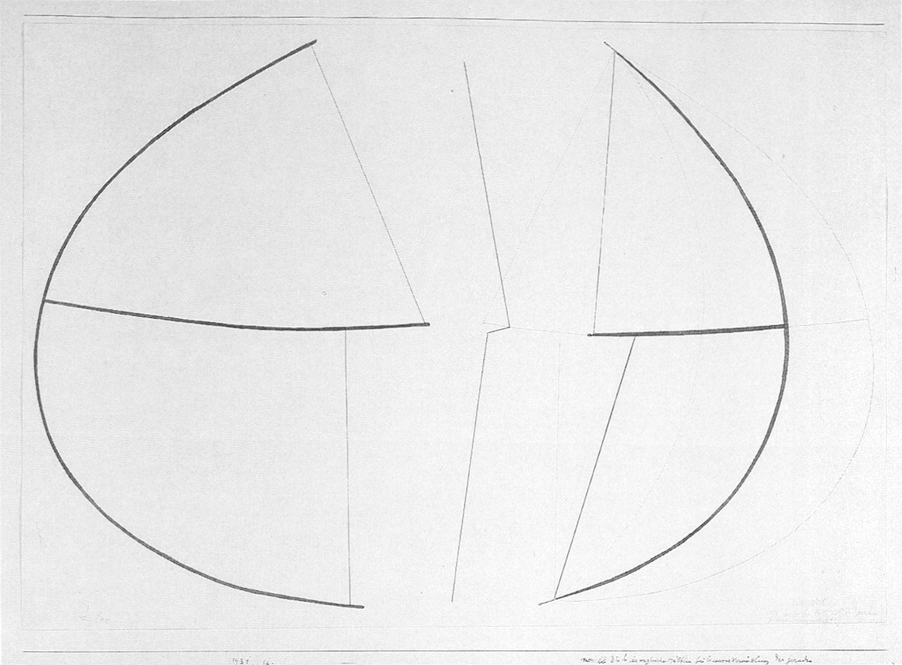}
\caption{\textit{Model 32B as Unequal Halves via Linear Mediation of Straight Lines} (1931), by Paul Klee, ink on paper mounted on board, 45 x 58 cm. Extended loan and promised gift of the Carl Djerassi Trust I.}
\label{Figure10}
\end{figure}

Any interpretation of the appearance of a shape (or an artwork for that matter), can be said to impose a certain structure or organization on the appearance of the shape. This structure is a result, a byproduct, of the parts that one perceives in the shape (or artwork)---it is an invention of the viewer. This highlights two things: (i) the same shape is a host of structures, each one invented from the parts that one chooses to see in it, and (ii) any structure imposed on a shape as a result of an interpretation is not permanent but changeable.

My main goal in this paper is to define and study topology on a shape based on interpretations of its appearance into parts. More specifically, how can the parts recognized in an aesthetic interpretation be used to define a topological structure for the shape in terms of them? Given that such a topological structure is possible, what are some of the mathematical properties that can be stated about this structure? What modifications or adaptations of classical topological concepts are necessary in order to accommodate shapes in their full pictorial and spatial potential? 

\begin{figure}[!ht]
\centering
\includegraphics[scale=0.8]{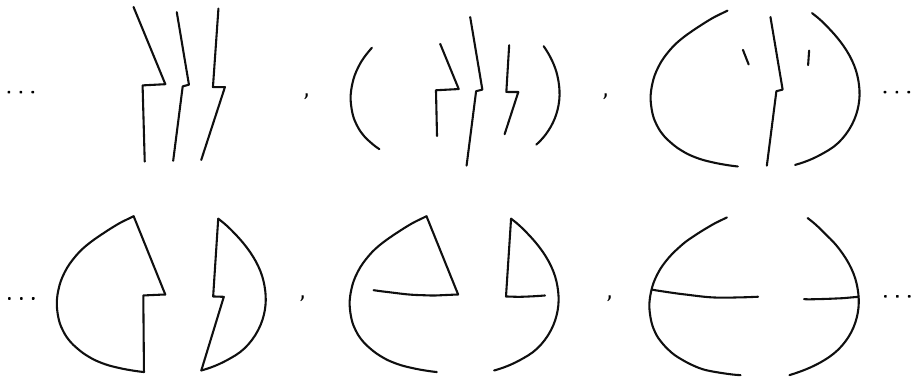}
\caption{Interpretations of Klee's artwork in Figure 10 terms of parts embedded in it.}
\label{Figure11}
\end{figure}

\begin{figure}[!ht]
\centering
\includegraphics[scale=0.25]{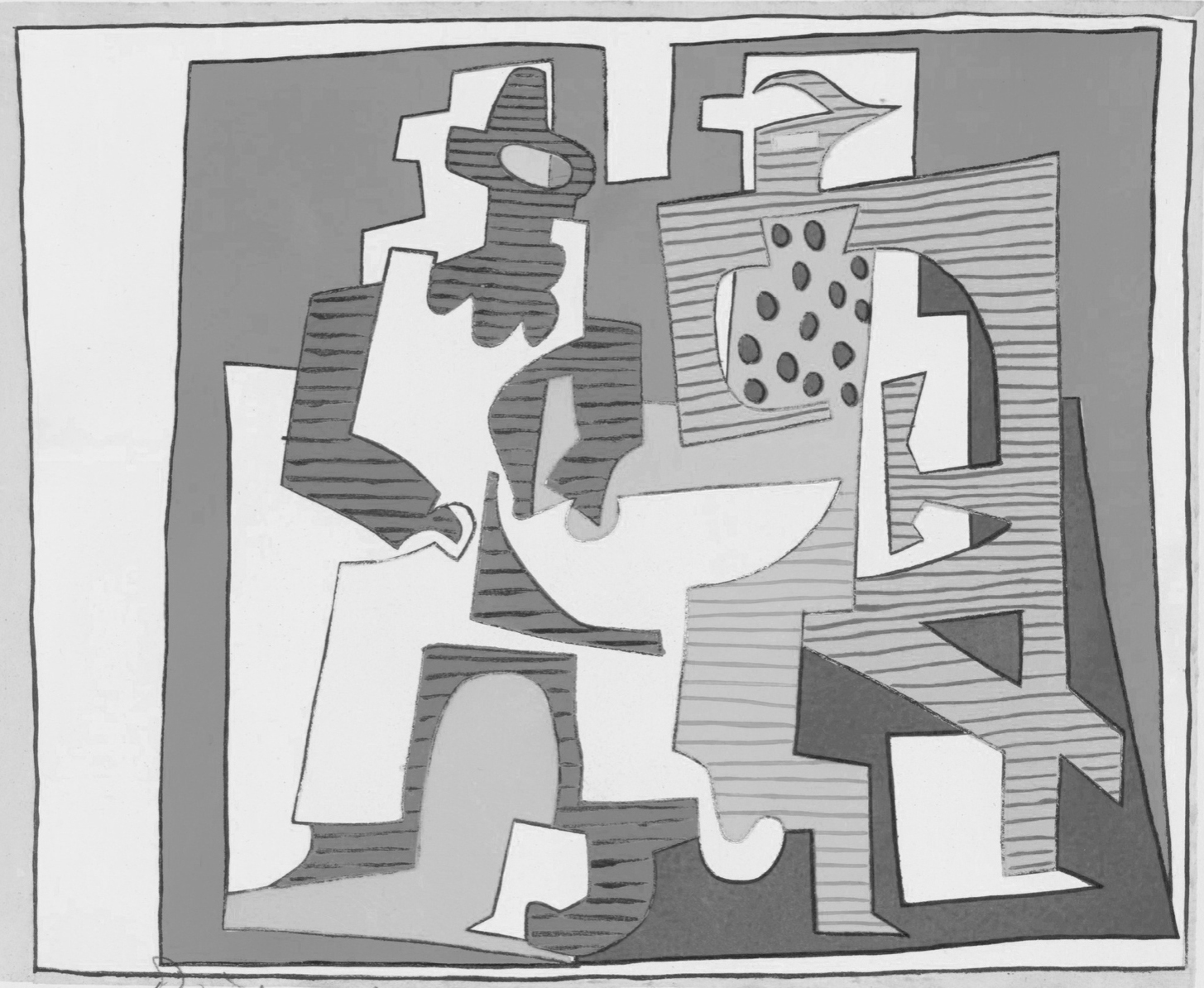}
\caption{\textit{Two Figures, Seated, after Pablo Picasso} (1920), by Pablo Picasso, stencil, 21.4 x 26.7 cm. \textcopyright\;2019 Estate of Pablo Picasso / Artists Rights Society (ARS), New York.}
\label{Figure12}
\end{figure}

\begin{figure}[!ht]
\centering
\includegraphics{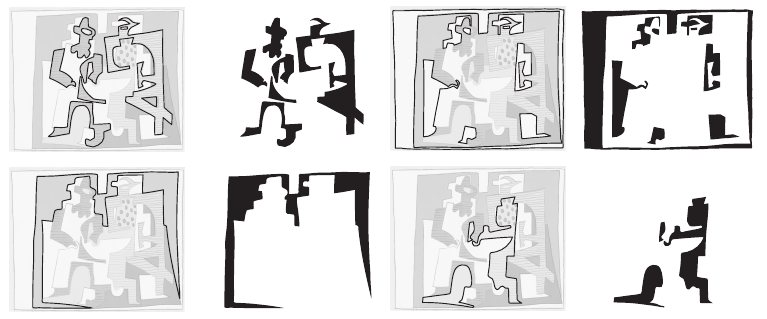}
\caption{Interpretations of Picasso's artwork in Figure 12 in terms of parts embedded in it. The parts are given as lines and as planes.}
\label{Figure13}
\end{figure}

\subsection{Previous work on finite topologies for shapes}

Topology for shapes in algebras $U_i$ is a relatively undeveloped research subject. In the literature of shape grammars, \cite{Stiny94} is the first study to propose a framework for topology by using closure algebras (historically, closure algebras---introduced in \cite{McKinTarski44}---were one of the first attempts in mathematics to define topology from algebraic information, without mentioning points \cite{Johnstone2001}). It is evident in this first study, that topologies for shapes can be defined as finite part-structures (partially-ordered structures that consist of finitely many parts). Mathematical properties of these structures were then studied in \cite{Krst96}. One of the main ideas in \cite{Krst96}, is that a topology for a shape can be described as a finite lattice of parts (the part relation in the algebras induces the order in the lattice).

A more general result about the topology underlying the algebras $U_i$ was given in \cite{Earl97}. For every algebra $U_i$, there is an associated topology given by the Stone representation of Boolean algebras (e.g., \cite{Johnstone82}). Any finite topology induced on a shape will essentially have its members (parts) embedded in the underlying Stone space of the specific algebra that the shape comes from. This underlying Stone space, however, adds no descriptive restrictions. That is to say, it does not signify, or restrict, the parts of a shape that can be distinguished for the purposes of a structural description \cite{Earl97, Stiny2006}.

Besides the aforementioned studies, no other significant work has been produced on the subject. This paper essentially provides the first exclusive and systematic treatment of finite topology for shapes in algebras $U_i$. The focus is on shapes made of basic elements of dimension $i > 0$. The topology applicable to shapes made with points (when $i = 0$) is not discussed here; it is the subject of a different paper \cite{Haridis2020}.

The main methodological approach is to replace the earlier formulation of finite topologies through closure algebras with a new formulation in which the concept of ``open part" is the only \emph{primitive} concept assumed. All topological concepts and constructions are built in terms of this assumption. Concepts that have been introduced in earlier publications are reconfigured (e.g., \emph{Sections 4} and \emph{5}), and many foundational concepts are introduced for the first time (e.g., \emph{Sections 3}, \emph{6}, and \emph{7}).

In the following presentation, topological ideas are formulated in such a way so that they can have an immediate and natural grounding in pictorial and spatial examples. Most of the time, the intuition comes from imagining what happens when graphical figures interact spatially. From a technical standpoint, this distinction is crucial because it is exactly what highlights the many ways in which shapes are different from point-sets and topological spaces.

\section{Shape topology}

\subsection{Definition and the concept of open part}

A \emph{topology} for a shape is a finite set $\mathcal{T}$ of parts of the shape that satisfies the following three conditions: 

\begin{enumerate}
\item[(1)] The \emph{empty shape} (\emph{0}) and the shape itself are in $\mathcal{T}$. 

\item[(2)] The sum (+) of an arbitrary number of parts in $\mathcal{T}$ is also in $\mathcal{T}$. 

\item[(3)] The product ($\cdot$) of an arbitrary number of parts in $\mathcal{T}$ is also in $\mathcal{T}$. 
\end{enumerate}

\noindent The members of a topology for a shape are shapes; they are the parts of the shape the topology recognizes. A shape without points has uncountably many parts (\emph{Section 2}). Thus, the requirement that a topology be finite implies that a shape without points can be induced with uncountably many different topologies, each one made up of \emph{finitely many} parts (infinite topologies on shapes are not considered in this paper). The set \{\emph{0}, $S$\} is the \emph{smallest} finite topology we can have on any given shape $S$. There exists no \emph{largest} finite topology on a shape without points.

Given a shape $S$, a topology can be constructed on $S$ based on the parts recognized in an interpretation of $S$. This motivates the following recursive description of a generated topology.

\vspace{12pt}

{\small \noindent \emph{Recursive description of a generated topology}\;\; Let $\mathcal{P}$ be a finite set of nonempty shapes all of which are part of shape $S$, and such that the sum of the shapes in $\mathcal{P}$ is equal to $S$. Define the sets $\mathcal{T}_0$, $\mathcal{T}_1$, $\mathcal{T}_2$... recursively as follows:

\begin{enumerate}
\item[(i)] $\mathcal{T}_0 = \mathcal{P} \cup \{ \emph{0} \}$. 

\item[(ii)] For each $k \geq 1$, define the set 

\vspace{0.1in}
\centerline{$\mathcal{T}_k$ = $\mathcal{T}_{k-1} \cup \{s_1, ..., s_m\}$}
\vspace{0.1in}

\noindent where $s_1, ..., s_m$ are shapes such that, for all $i = 1, ..., m$, shape $s_i$ is either a sum of a finite number of nonempty shapes in $\mathcal{T}_{k-1}$, or a product of such shapes.

\end{enumerate}

\noindent The topology on $S$ generated by $\mathcal{P}$ is equal to $\mathcal{T}_{n}$ for a $n < \infty$ for which $\mathcal{T}_{n} = \mathcal{T}_{n+1}$, that is to say $\mathcal{T}_{n+1}$ introduces no new parts.}

\vspace{12pt}

This recursive description provides a procedure for obtaining a topology in terms of certain initial parts, in finitely many steps ($\mathcal{P}$ is finite and (i) and (ii) are determined by finite operations). A concrete example is given in \emph{Example 1}. 

The same procedure can be followed even if the shape $S$ is already equipped with a topology $\mathcal{T}$. Suppose $\mathcal{P}$ is a set of newly recognized parts of $S$, none of which is already open in the existing topology for $S$. Then, let $\mathcal{T}_0$ be equal to $\mathcal{T} \cup \mathcal{P}$ and proceed recursively as described above (it is assumed that the \emph{empty shape} is already included in $\mathcal{T}$). The resulting topology will be a \emph{refinement} of the existing topology on $S$ in terms of the newly recognized parts. 

A few notes on the terminology. At times I refer to the topology $\mathcal{T}$ for a shape $S$ as \emph{shape topology} (this is for distinguishing between topologies for shapes and topologies for parts, i.e., subshape topologies; see \emph{Section 3.3}). To distinguish between the parts of a shape that are members of a topology from those that are not, I use the term \emph{open parts} to refer to the parts in the topology. Thus, as a general rule of thumb, on the side of the terminology, the following two-way implication is assumed:

\vspace{12pt}

\centerline{$x$ is a recognized part of $S$ $\iff$ $x$ is \emph{open} in $S$, i.e., $x$ is in a topology for $S$.}

\vspace{12pt}

\noindent Using this terminology, one can say that a topology for a shape $S$ is a set $\mathcal{T}$ of open parts of $S$, satisfying the three basic conditions for a topology. 

The notation ($S$, $\mathcal{T}$) can be used to refer to a shape $S$ that possesses a topology $\mathcal{T}$. Alternatively, the notation $OS$ can be used to refer specifically to the \emph{lattice of open parts} determined by the topology $\mathcal{T}$ (the order of this lattice is induced by the part relation ($\leq$) in the algebras $U_i$, where $i > 0$; this lattice is complete, distributive, with join and meet replaced, respectively, by sum and product for shapes). Terminological conventions brought by lattice theory are in general avoided, except for the cases where those conventions are important; see, for example, \emph{Section 6} on \emph{Continuity}.

The parts of a shape which are members of a particular topology are the \emph{only} open parts of this shape with respect to that topology. Any part of a shape is ``visible", in the sense that we can see it, but is not necessarily open. A part of a shape is open, if and only if, it is a member of a topology for the shape.

Finite topologies for shapes can be represented in two ways. First, by a list of open parts (with an implied ordering). Second, by a lattice diagram. In this paper, I use the latter representational method throughout. To simplify the presentation, I have chosen to omit directed lines between open parts.

\subsection{Basis for a topology on a shape}

It is possible to specify a topology for a shape by describing the entire set $\mathcal{T}$ of open parts. When this is too difficult, one specifies a smaller set of parts and describes the topology in terms of them. This smaller set is called a \emph{basis}.

Let $S$ be a shape and $\mathcal{B}$ a set of parts of $S$. $\mathcal{B}$ is a basis for a topology on $S$ if:

\begin{enumerate}
\item[(1)] The sum of the parts in $\mathcal{B}$ is equal to $S$.

\item[(2)] If $b_1$ and $b_2$ are two parts in $\mathcal{B}$, then there is a third part $b_3$ in $\mathcal{B}$, such that $b_1 \cdot b_2$ = $b_3$, or $b_1 \cdot b_2$ is the sum of parts from $\mathcal{B}$.
\end{enumerate}

\noindent If $\mathcal{B}$ satisfies conditions (1) and (2), we define the topology $\mathcal{T}$ generated by $\mathcal{B}$ as follows: A part $x$ of $S$ is said to be open (that is, to be a member of $\mathcal{T}$) if it is a basis element of $\mathcal{B}$ or if it can be described as a sum of basis elements from $\mathcal{B}$. 

From this defining condition of openness it is obvious that all basis elements are themselves open. Moreover, given a part $C$ open in $S$, if $C$ is not itself a basis element, then we can choose basis elements $b_i$ embedded in $C$, i.e., $b_i \leq C$, so that $C = \sum b_i$; this fact proves very handy in many occasions in this and in the following sections. Note the \emph{empty shape} is always included as a member of a basis. It will be usually omitted, however, whenever a basis is presented explicitly. 

Let us now check that the set $\mathcal{T}$ generated by $\mathcal{B}$ is indeed a topology for $S$ (i.e., that it satisfies the three conditions given in \emph{Section 3.1}). The shape $S$ is open, by condition (1) above. Similarly, the \emph{empty shape} is open because it is a member of $\mathcal{B}$, by definition. Now, given $C = \sum C_\alpha$ an arbitrary sum of open parts, $C$ is open because for each index $\alpha$, $C_\alpha$ is either a basis element or the sum of basis elements, by definition. 

Finally, we need to show that $D = \prod D_\alpha$, an arbitrary (finite, in our case) product of open parts, is open. Let us illustrate this in the simpler case where $D$ = $D_1 \cdot D_2$. Since $D_1$ and $D_2$ are open, they are either basis elements or can be expressed as sums of basis elements (according to the defining condition of openness). Without loss of generality, suppose $D_1 = \sum b_i$ and $D_2 = \sum b_j$. Then, we can write $D$ as

\vspace{12pt}
\centerline{$\sum b_i \cdot \sum b_j = \sum (b_i \cdot b_j)$,} 
\vspace{12pt}

\noindent for all $i, j$. By condition (2) above, the product $b_i \cdot b_j$, for all $i, j$, is either a basis element or a sum of basis elements, so that $D$ is open as desired.

We now know that to describe a topology for a shape $S$, it suffices to present a set $\mathcal{B}$ of parts of $S$ which satisfies conditions (1) and (2). To actually construct the topology, that is to say, to compute all the open parts, the recursive procedure of \emph{Section 3.1} can be used; simply let $\mathcal{P}$ be equal to $\mathcal{B}$ and proceed in the obvious way (from conditions (1) and (2), it is easy to see that sum is the only operation needed in this case).

Often a topology for a shape can be described with more than one set of parts. That is to say, a topology may have multiple bases (for example, even the topology itself can be taken as a basis). Since the topologies we deal with here are finite, we can naturally describe them with ``minimal" bases. Let's examine what is exactly expected from such a basis. 

A minimal basis is a basis with the \emph{smallest} possible cardinality or, equivalently, a basis whose elements are contained in every other basis for the same topology (obviously, such a basis must be unique for finite structures). In principle, it must be able to describe (i.e., generate) the elements of every other basis for the same topology. This can be achieved by forming a set of parts that does \emph{not} contain a part which can be expressed as a sum of other parts from the \emph{same} set. It is not hard to see, if not immediate, that a basis that satisfies conditions (1) and (2) only, is not necessarily minimal---we can often create multiple bases for a topology, all of which satisfy the two conditions, but none of which is minimal. Even if we are presented with a minimal basis for a particular topology, we currently have no general rule for telling if it is actually minimal or not (unless of course if we compare it with every other possible basis, one by one).

To address minimality, an additional requirement is introduced. A basis $\mathcal{B}$ for a topology on $S$ is a \emph{minimal basis} if, in addition to (1) and (2), it also satisfies the following condition:

\begin{enumerate}
\item[(3)] If $b_1$ and $b_2$ are two parts in $\mathcal{B}$ and $b_1$ + $b_2$ is also in $\mathcal{B}$ then either $b_1$ + $b_2$ = $b_1$ or $b_1$ + $b_2$ = $b_2$.
\end{enumerate}

\noindent In \emph{Appendix A.1}, I show that for any topology $\mathcal{T}$, a set of parts $\mathcal{B}$ that satisfies conditions (1) through (3) is a \emph{unique minimal basis} for $\mathcal{T}$. If we are presented with a basis for a topology we can always reduce it to its (unique) minimal version by removing the elements that do not obey condition (3). Moreover, we can use these conditions to decide if a given basis is minimal or not (e.g., see \emph{Example 1}).

In the rest of this paper, I will refer to the unique minimal basis of a topology as a \emph{reduced basis}; whether a basis is supposed to be reduced or not will be made explicit every time from context.

\subsection{Subshape topology and covering}

Suppose $S$ is a shape with a topology $\mathcal{T}$. This topology can be ``transferred", or relativized, to a part of $S$ that is not necessarily open in $\mathcal{T}$. This part then inherits a topology determined by the existing open parts in $\mathcal{T}$. This motivates the concept of a \emph{subshape topology}. Let $x$ be a part of $S$. The set,
\[
\mathcal{T}_x = \{\;x \cdot C \;|\; C \textrm{ an \emph{open part} in } \mathcal{T}\;\}
\]
\noindent is a topology for $x$, called the subshape topology. The open parts in $\mathcal{T}_x$ consist of all products of open parts of $S$ with the part $x$. 

Let us check that the set $\mathcal{T}_x$ is a topology. The shape $x$ and the \emph{empty shape} are in $\mathcal{T}_x$ because

\vspace{0.1in}
\centerline{\emph{0} = $x \cdot \emph{0}$ and $x$ = $S \cdot x$,}
\vspace{0.1in}

\noindent where \emph{0} and $S$ are open in $\mathcal{T}$. The fact that arbitrary sums and products of open parts are open in $\mathcal{T}_x$ follows from the equations 

\vspace{0.1in}
\centerline{$\sum (C_a \cdot x) = (\sum C_a) \cdot x$ and $\prod (C_a \cdot x) = (\prod C_a) \cdot x$.}
\vspace{0.1in}

Using subshape topology, one can pick an arbitrary part $x$ of $S$ and construct a topology on $x$ by reusing the existing topology $\mathcal{T}$ on $S$. The subshape topology is said to be a relativization of $\mathcal{T}$ with respect to the part $x$. If $x$ is already open in $\mathcal{T}$, then obviously $\mathcal{T}_x \subset \mathcal{T}$.

Similar to the basis for a shape topology, we can define a basis for the subshape topology. If $\mathcal{B}$ is a basis for a topology $\mathcal{T}$ on $S$, then the set, 
\[
\mathcal{B}_x = \{\;x \cdot b \;|\; b \textrm{ a \emph{basis element} in } \mathcal{B}\;\}
\]
\noindent is a basis for the subshape topology $\mathcal{T}_x$ on $x$; the members of $\mathcal{B}_x$ consist of all products between $x$ and all the basis elements in $\mathcal{B}$. That the set $\mathcal{B}_x$ is indeed a basis for $\mathcal{T}_x$, is proven in \emph{Appendix A.2}. In \emph{Example 1}, I show how to construct subshape topologies based on an existing topology for a shape.

When referring to open parts, one needs to specify if those parts are open with respect to a subshape topology or with respect to a shape topology. By transitivity, if a part is open in the subshape topology $\mathcal{T}_x$ for a part $x \leq S$ and $x$ is itself open in a topology $\mathcal{T}$ for $S$, then this part is also open in $\mathcal{T}$. Unless otherwise stated, an open part will be always assumed to be ``open" with respect to a (shape) topology for $S$.

Some further topological concepts are readily defined. A set of open parts of $S$ is said to \emph{cover} $S$, or to be a \emph{covering} of $S$, if the sum of the parts in this set is equal to $S$. For example, a basis for a topology on $S$ is automatically a covering for $S$. A \emph{subcovering} of $S$, is a subset of a covering that still covers $S$. A covering of a nonempty part $x$ of $S$ is a set of open parts of $S$, such that the sum of the parts in this set has $x$ as a part. A subcovering for $x$ is defined analogously.

\begin{figure}[ht]
\centering
\includegraphics{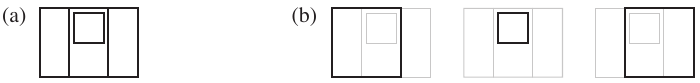}
\caption{(a) A shape and (b) three selected visible parts of this shape.}
\label{Figure14}
\end{figure}

A complementary idea is the intuitive notion of ``exhaustion". Suppose you want to recognize certain parts in a given shape $S$. If the sum of these parts is equal to $S$, then these parts can be said to exhaust the (appearance of) shape $S$. For example, the parts recognized in Figure 14b exhaust the shape in Figure 14a. In \emph{Section 2}, the parts recognized in all rows of Figure 8 exhaust the shape in Figure 6a. It is not, however, necessary for an interpretation to be always exhaustive; one may choose to recognize only some portion of $S$ (not the whole shape). For example, the right shape in Figure 6c, and the series of shapes in Figure 9, do not exhaust the shape in Figure 6a. In that case, the recognized parts do not sum to $S$, and hence, do not exhaust $S$. This concept of exhaustion is now connected to topology.

A topology $\mathcal{T}$ for a shape $S$ is said to \emph{exhaust} $S$, when the sum of the parts in the set $\mathcal{T} \setminus \{S\}$ is equal to $S$. (The symbol $\setminus$ stands for the operation of set-difference.) If $\mathcal{T}$ does not exhaust $S$, then the sum of the parts in the set $\mathcal{T} \setminus \{S\}$ leaves an undivided complement relative to $S$. For example, the topologies in Figure 15a, Figure 16a and 16b exhaust the shape in Figure 14a; the topology in Figure 16c does not. 

When a topology $\mathcal{T}$ does \emph{not} exhaust $S$, the basis that generates $\mathcal{T}$ must include $S$ itself. In other words, if a topology does not exhaust a shape, the shape itself has to be taken as one of its basis elements.

\begin{figure}[!ht]
\centering
\includegraphics[width=\textwidth]{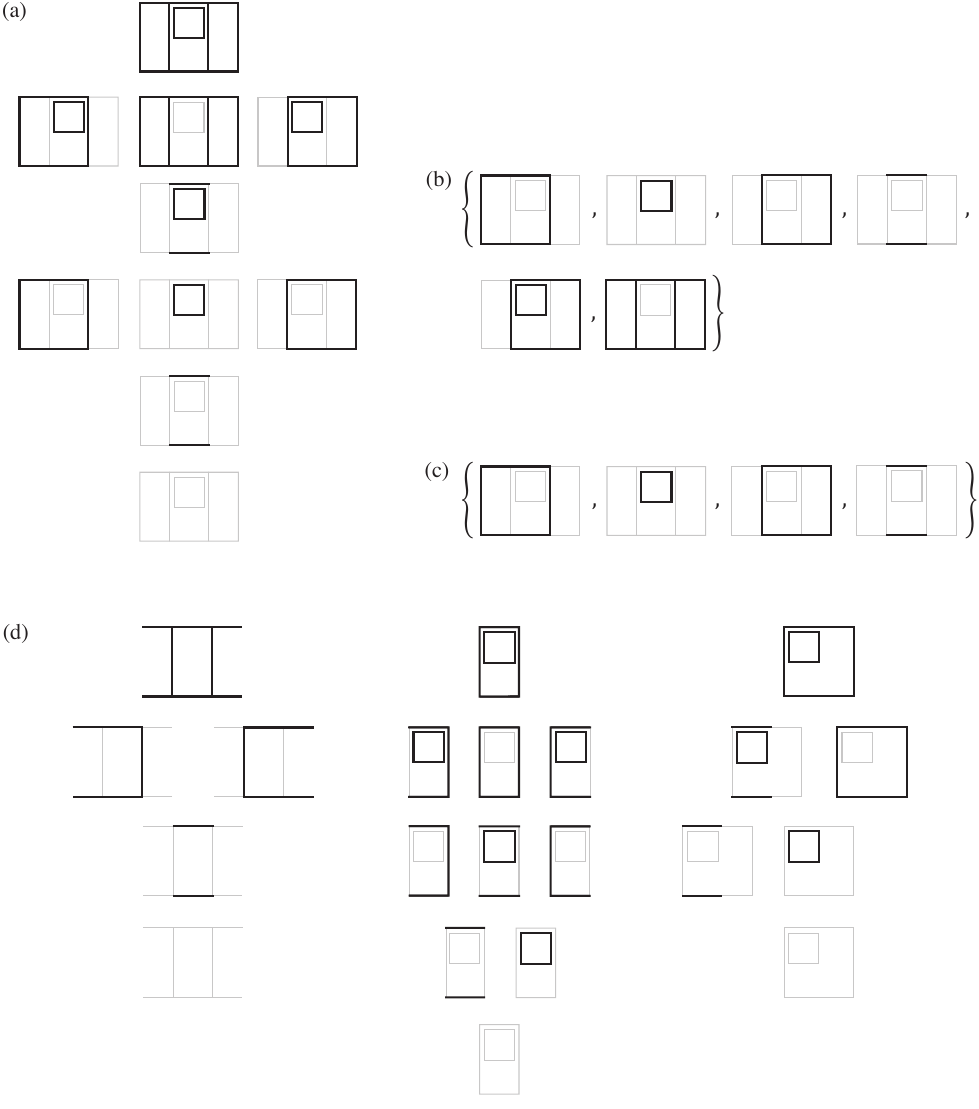}
\caption{(a) A topology induced on the shape in Figure 14a. (b) and (c) are two bases that generate this topology, where (c) is in reduced form. (d) Three subshape topologies of the shape topology in (a).}
\label{Figure15}
\end{figure}

\vspace{12pt}
\noindent \emph{EXAMPLE 1. \;\;\;Let $S$ be the shape in Figure 14a. Suppose we choose to recognize in $S$ the three parts shown in Figure 14b. A topology $\mathcal{T}$ can be defined on $S$ in terms of these three parts using the recursive procedure given in Section 3.1. Figure 15a shows the resulting topology on $\mathcal{T}$. The sets of parts in Figure 15b and 15c both generate $\mathcal{T}$---thus, both sets are bases for $\mathcal{T}$. Only the set in Figure 15c, however, satisfies the three conditions for a reduced basis for $\mathcal{T}$. Figure 15d shows three different subshape topologies for three different parts of the shape $S$. The reader may want to verify that the structures in Figure 15d are indeed subshape topologies; the computation of the reduced bases that generate each of them is left as an exercise.}



\subsection{Comparing shape topologies}

Two or more topologies induced on the same shape can be compared with respect to the sets of parts that each one recognizes. Let $\mathcal{T}_1$ and $\mathcal{T}_2$ be two topologies induced on the same shape. Say that $\mathcal{T}_2$ is \emph{finer} (or larger) than $\mathcal{T}_1$ if every part in $\mathcal{T}_1$ is also in $\mathcal{T}_2$, or in notation, $\mathcal{T}_2 \supset \mathcal{T}_1$. In this case, $\mathcal{T}_1$ is said to be \emph{coarser} than $\mathcal{T}_2$. Two topologies $\mathcal{T}_1$ and $\mathcal{T}_2$ are said to be \emph{comparable} if either $\mathcal{T}_2 \supset \mathcal{T}_1$ or $\mathcal{T}_2 \subset \mathcal{T}_1$ holds; they are not comparable if neither inclusion holds. When both inclusions hold at the same time, $\mathcal{T}_1 = \mathcal{T}_2$ in which case $\mathcal{T}_1$ and $\mathcal{T}_2$ recognize exactly the same parts. Figure 16 shows three different topologies for the same shape. Topology $\mathcal{T}_2$ (Figure 16b) is finer than $\mathcal{T}_1$ (Figure 16a), but $\mathcal{T}_1$ and $\mathcal{T}_2$ are not comparable to $\mathcal{T}_3$ (Figure 16c).

\begin{figure}[ht]
\centering
\includegraphics[width=\textwidth]{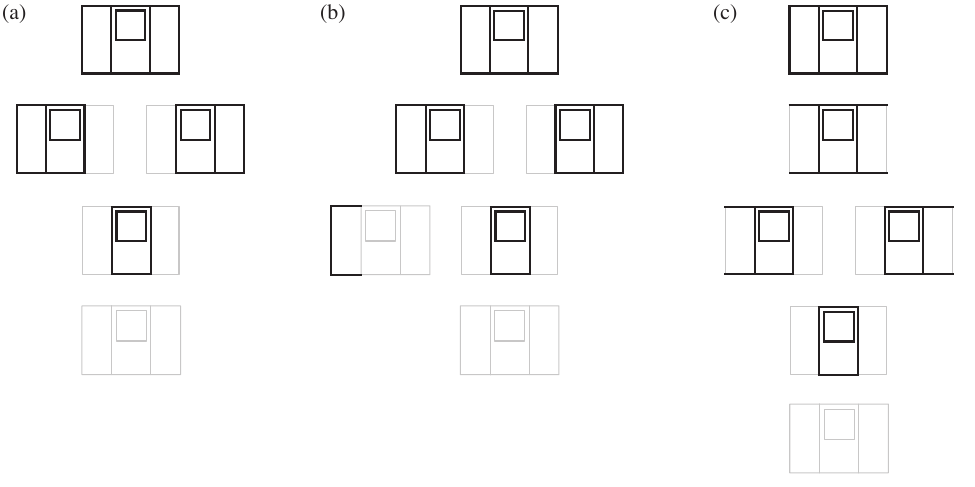}
\caption{Three different shape topologies induced on the shape in Figure 14a. The topologies (a) and (b) are comparable; topology (c) is not comparable to either (a) or (b).}
\label{Figure16}
\end{figure}

Comparisons of topologies are in a sense comparisons of ``granularities". A comparison of two topologies is often straightforward to make: simply count the number of open parts in the two given topologies and check if the topology with the fewer pieces is entirely included in the topology with the most pieces. When the topologies are instead described by their reduced bases, the following criterion can be used to compare them in terms of those bases. 

Let $\mathcal{B}_1$ and $\mathcal{B}_2$ be the reduced bases that generate topologies $\mathcal{T}_1$ and $\mathcal{T}_2$, respectively, on the same shape. Then, the following are equivalent:

\begin{enumerate}
\item[(1)] $\mathcal{T}_2$ is finer than $\mathcal{T}_1$.

\item[(2)] For every basis element $b$ in $\mathcal{B}_1$, there is a nonempty subset of basis elements in $\mathcal{B}_2$ such that the sum of the elements in this subset is equal to $b$.
\end{enumerate}

\noindent For example, criterion (2) gives the correct comparability result for topologies $\mathcal{T}_1$ and $\mathcal{T}_2$ in Figure 16, i.e., that $\mathcal{T}_2$ is finer than $\mathcal{T}_1$; for every element in the (reduced) basis of $\mathcal{T}_1$ one can find a subset of elements in the (reduced) basis of $\mathcal{T}_2$ that sums to this element, but the opposite is not possible. Criterion (2) also gives the correct compatibility results between $\mathcal{T}_1$ and $\mathcal{T}_3$, and $\mathcal{T}_2$ and $\mathcal{T}_3$ in the same figure. Namely, that $\mathcal{T}_3$ is not comparable to either $\mathcal{T}_1$ or $\mathcal{T}_2$.

Comparisons of finite topologies on shapes are interesting from an interpretative, aesthetic standpoint. A shape can be interpreted aesthetically in a variety of ways, depending on the parts that different viewers choose to see in it. The different topologies derived from those parts do not have to be comparable with one another---indeed, the parts seen by a viewer in a shape do not have to ``match" with the parts seen by another viewer in the same shape (or even by the same viewer at a different time). Thus, a shape can be a host of non comparable, possibly conflicting and contradictory structures.

\section{Closed part, interior and closure}

We have defined what topology is for a shape without points, how to generate topologies with bases, and how to compare different topologies. Now we proceed to concepts that are useful for describing parts of a shape based on their topological characteristics.

Let $S$ be a shape with topology $\mathcal{T}$. A part $x$ of $S$ is said to be \emph{closed} whenever its complement relative to $S$, namely $S$ - $x$, is open. For example, in Figure 17, part (a) is closed but not open in the topology in Figure 15a; likewise, parts (b), (c) and (d) are closed but not open in the first, second and third topology in Figure 15d, respectively.

\begin{figure}[ht]
\centering
\includegraphics{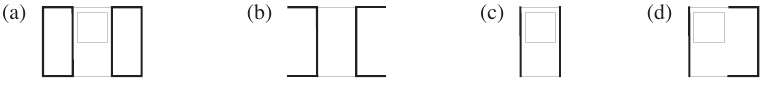}
\caption{Parts that are closed but not open in, respectively, (a) the topology in Figure 5a, (b) the first, (c) second and (d) third topologies in Figure 5d.}
\label{Figure17}
\end{figure}

Instead of using the notion of ``open part" as the primitive concept, one could just as well define a topology for a shape by using the notion of ``closed part". One could then define open parts as the complements of closed parts and proceed in the obvious way. This procedure has no particular technical advantages over the one already adopted.

A part $x$ of $S$ is said to be \emph{closed-open} if $x$ is both open and closed in $\mathcal{T}$. In other words, a closed-open part is an open part whose complement is also open. For example, in any topology the \emph{empty shape} and the shape itself are always closed-open parts. The topology in Figure 15a, the second and third topologies in Figure 15d have additional closed-open parts---they are shown in Figure 18a, 18b and 18c, respectively. The first topology in Figure 15d has no additional closed-open parts.

\begin{figure}[ht]
\centering
\includegraphics{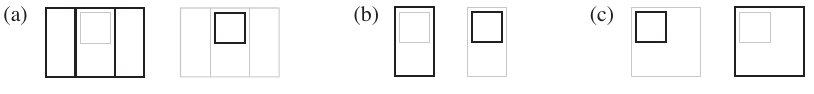}
\caption{Closed-open parts in, respectively, (a) the topology in Figure 15a, (b) the second and (c) third topologies in Figure 15d.}
\label{Figure18}
\end{figure}

The \emph{interior} of a part $x$ of $S$ is the sum of all the open parts in $\mathcal{T}$ that are embedded in $x$. One easily notices that the interior of $x$ is the \emph{largest member} in $\mathcal{T}$ that is embedded in $x$. It follows that if $x$ is open, the interior of $x$ is equal to $x$ itself. For example, in any topology, $S$ is always equal to its interior. If $x$ is not open, the interior of $x$ is a proper, possibly empty, part of $x$.

The following is a particularly useful result concerning finite topologies for shapes. For any part $x$ of $S$, there is a unique \emph{smallest member} in $\mathcal{T}$ that has $x$ as a part. 

\begin{proof}
Let $C$ be the product of all open parts of $S$ in $\mathcal{T}$ which have $x$ as a part. Since $\mathcal{T}$ is finite, this is a finite product and so $C$ is open. It is immediate that $C$ must be the smallest member in $\mathcal{T}$ that has $x$ as a part and that $C$ is unique.\end{proof}

The existence of the unique smallest member supports the following definition, which is dual to the concept of interior. The \emph{closure} of a part $x$ of $S$ is the smallest member in $\mathcal{T}$ that has $x$ as a part. For any part, its closure is determined each time by the available open parts in a topology. The following is obtained as a consequence of this: a part $x$ of $S$ is open if and only if $x$ is equal to its closure in $\mathcal{T}$.

To introduce some notation, let $Int x$ denote the interior of a part $x$ and $\overline{x}$ its closure. The following chained relationship describes how the interior of the part, the part itself, and the closure of the part are related:
\[
Int x \leq x \leq \overline{x}
\]
The following figure is a pictorial illustration of this relationship for an arbitrary part of the shape in Figure 14a; its interior and closure are determined based on the topology in Figure 15a.

\begin{figure}[!ht]
\centering
\includegraphics{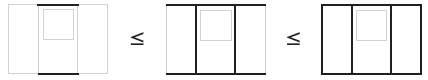}
\label{interior-closure}
\end{figure}

One may expect to name a part that is equal to its closure a ``closed part" (for example, this is how closed parts are defined in \cite{Stiny94}). However, this is not in sync with how finite topologies are developed in this paper. The name closed part is reserved for a part whose \emph{complement} is open (but is not necessarily itself open). Additionally note that \emph{any} visible part of a shape, arbitrarily chosen, can have a closure in a topology, independently of whether this part, or its complement, is open or not.

It is clear from the given definitions of interior and closure that if a part $x$ of $S$ is open, then $x$ is equal to both its interior and its closure. This indicates that there is no real distinction between the interior of a part, the part itself, and its ``limits", as there is, for example, for open sets and their limit points in a topological space. For shapes, the closest concept to the notion of limit is perhaps the boundary. Indeed, shapes do have boundaries. However, as explained in \emph{Section 2.1}, the boundary of a shape is not made from the same basic elements as the shape itself, and thus it is not part of the shape. The closure of a part that is open ends up simply being the part itself; if a part is not open, its closure is the smallest open part in the topology that this part is embedded in.

\section{The \emph{space} of a shape relative to a topology}

The set of points onto which one induces a topological structure is known in general topology as \emph{space}. A space can be infinite, such as the real number line $\mathbf{R}$ or the real plane $\mathbf{R}^2$, or finite, such as the set of labelled vertices of a finite graph. No matter the choice of the space, it is assumed that the points are given all at once, ahead of time. Moreover, the points of the space are absolute, in the sense that they are invariant under different topologies. For example, regardless if the real plane $\mathbf{R}^2$ is equipped with a metric topology or with the standard topology generated by open rectangles, the points of $\mathbf{R}^2$ are the same. Is this also true for shapes and shape topologies?

The first publication that considered the notion of ``space" in conjunction with finite topology on shapes in algebras $U_i$ is \cite{Krst96}. As it is rightly pointed out there, for any given shape, an underlying space can be determined only \emph{relative to} a specific topology. In other words, a topology induced on a shape gives rise to its own underlying space and thus different topologies for the same shape will give rise to different spaces. To actually define the space of a shape relative to a topology, \cite[p. 128]{Krst96} suggests mathematical machinery originating in the theory of locales (point-free topology, e.g., \cite{PicPultr2010}). Here, an equivalent, albeit much simpler, approach is suggested, using only material developed in \emph{Section 3} and without appeal to extra machinery: If $S$ is a shape and $\mathcal{T}$ a topology for $S$, the \emph{space} underlying $S$ relative to $\mathcal{T}$ is precisely the reduced basis $\mathcal{B}$ that generates $\mathcal{T}$, and the basis elements of $\mathcal{B}$ may be called the \emph{points} of this space (i.e., the ``points" of $S$ relative to $\mathcal{T}$; the points obtained from this approach correspond to the f-prime elements of a topology in the formalization of topologies for shapes within locale theory that \cite{Krst96} follows).

We can cast any topology $\mathcal{T}$ for a shape into its space theoretical version, say $\mathcal{T}^*$, using the reduced basis for $\mathcal{T}$. This is easily shown in an example.

\begin{figure}[ht]
\centering
\includegraphics[width=\textwidth]{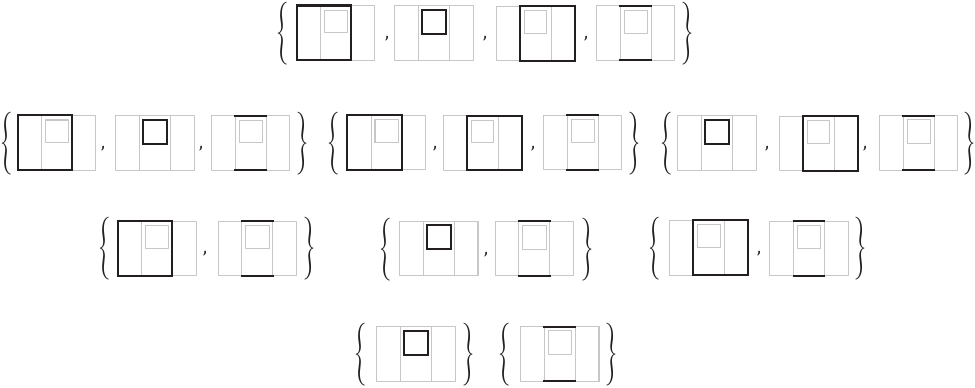}
\caption{A topology induced on the space for the shape in Figure 14a. The space and the topology are constructed relative to the topology in Figure 15a. Each member of this topology is a set and corresponds to an open part from the original topology.}
\label{Figure19}
\end{figure}

Suppose $S$ is the shape in Figure 14a and $\mathcal{T}$ the topology for $S$ in Figure 15a. The space of $S$ relative to $\mathcal{T}$ is the reduced basis $\mathcal{B}$ in Figure 15b. Let $C$ be any nonempty open part in $\mathcal{T}$. Define the set of elements of $\mathcal{B}$ embedded in $C$,
\[
B_C = \{\; b \in \mathcal{B} \;|\; b \leq C\;\}.
\]
Then, the collection $\{B_C\}_{C \in \mathcal{T}}$, with the empty set ($\emptyset$) included, forms a topology $\mathcal{T}^*$ on the space of $S$. This topology is shown graphically in Figure 19. Notice that each member $B_C$ of $\mathcal{T}^*$ is no longer a shape but a set; the elements in each set sum to the open part $C$ coming from $\mathcal{T}$. Thus, if we compute the sum of the elements in all sets (and map the empty set to the \emph{empty shape}), the resulting set of shapes will be equal to the original topology $\mathcal{T}$. By construction, $\mathcal{T}$ and $\mathcal{T}^*$ are isomorphic structures.

For any shape, there is no single absolute space underlying all finite topologies that can be induced on it. Rather, every different topology induced on a shape specifies its own underlying space via its reduced basis. In contrast, for shapes made with points, the points of the shape stay absolute and independent of the topologies induced on it, just like in topology for finite sets \cite{Haridis2020}.

\section{Continuity}

\subsection{Mappings between shapes}

Continuity is an important concept and one of the foundational ones in topology. To study continuity, one needs some appropriate notion of mapping (or function) between the objects considered.

\begin{figure}[ht]
\centering
\includegraphics{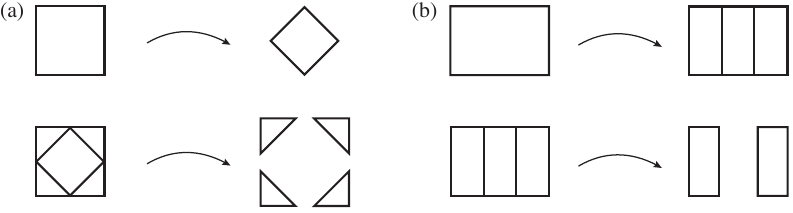}
\caption{Pictorial illustration of mappings between shapes.}
\label{Figure20}
\end{figure}

A mapping between two shapes can be thought of as a certain ``rule of assignment" that assigns/maps \emph{every part} of one shape to a part of the other shape. More intuitively, a mapping is a description of a \emph{spatial transformation} of the parts of one shape into parts of another shape. To emphasize this spatial aspect of mappings we represent them with graphical diagrams, such as those in Figure 20. 

A mapping may describe how a shape (the whole shape), or part(s) of it, is transformed by rotation, reflection, translation, or a linear transformation (or other known transformations) into another shape; for example, the mappings in Figure 20a. A mapping may also describe how a shape is transformed to another by adding something to it (Figure 20b, top), or by subtracting (erasing) something from it (Figure 20b, bottom). A mapping may also take a shape into itself untouched, in which case we have a special kind of mapping known as the \emph{identity} mapping (Figure 21b). Mappings between shapes are used all the time in practice in art and design. They are one of the most popular ways of creating something new out of something already known.

A \emph{mapping} from a source shape $S$ into a target shape $S'$ is written as

\vspace{12pt}
\centerline{ $f: S \rightarrow S'$.}
\vspace{12pt}

\noindent If $x$ is a part of $S$, denote by $f(x)$ the (unique) image of $x$ determined under $f$; the image of $x$ is a part $y$ of $S'$, such that $f(x) = y$. The mapping $f$ describes how the parts of $S$ are spatially changed or transformed into parts of the target shape $S'$, so that the image of $S$ under $f$ is equal to the full shape $S'$, i.e., $f(S) = S'$. Because the set of all the parts of a shape forms a partially ordered set (poset) (see \emph{Section 2}), a mapping between two shapes can be understood as a mapping between their respective posets (the mappings we're interested in are only those that preserve the embedding order of parts, as explained in a forthcoming paragraph).

New mappings between shapes can be formed from given ones. For example, a new mapping can be formed by restricting the target image of $f$ to a \emph{proper part} of $S'$. In particular, the mapping $h: S \rightarrow S^{+}$, where $S^{+}$ is a proper part of $S'$ and $h(S) = S^{+}$, describes how $S$ is transformed to some part $S^{+}$ of $S'$, but not to the full shape $S'$. The \emph{composite} of two given mappings can be also formed in the obvious way: given $f: S \rightarrow S'$ and $g: S' \rightarrow S''$, the composite of the two is the mapping $(g \circ f): S \rightarrow S''$ defined by the equation $(g \circ f)(x) = g(f(x))$, for all parts $x$ of $S$. 

If $y$ is a part of $S'$, denote by $f^{-1}(y)$ the \emph{inverse image} or \emph{preimage} of $y$. Note, here the notation $f^{-1}$ represents an \emph{operation}, namely, the operation of preimage. It should not be confused with the inverse mapping of $f$ (i.e., the mapping $f^{-1}: S' \rightarrow S$). 

In set theory, we have that the image of a set (or of a subset) is a set, and the preimage of set (or of a subset) is also a set. A desired scenario for shapes would be: the image of a shape (or of a part) is a shape (as we already said), and the preimage of a shape (or of a part) is also a shape. It turns out that preimages of shapes don't behave as expected. Extra care is needed to obtain a definition suitable for our purposes. The definition we seek can be grounded nicely in a simple pictorial example.

Suppose that the following mapping between two shapes $S$ and $S'$

\begin{figure}[!ht]
\centering
\includegraphics{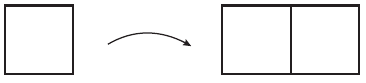}
\label{map-Example}
\end{figure}

\noindent is described by a mapping $f: S \rightarrow S'$, defined by $f(x) = x + A$, for any part $x$ of $S$, so that $f(S) = S'$. The symbol $A$ represents the shape 

\begin{figure}[!ht]
\centering
\includegraphics{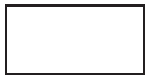}
\label{shape-A}
\end{figure}

\noindent that is added to $S$ to produce $S'$. Then, suppose we ask for the preimage of this part $y$ of $S'$,

\begin{figure}[!ht]
\centering
\includegraphics{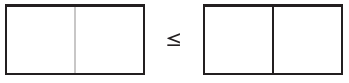}
\label{mapPart-Y}
\end{figure}

\noindent Following an approach analogous to the classical definition of preimages in set theory (e.g., \cite{Munkr00}), we need to define a set of parts of shape $S$ all of which have the following property: the image of each part in this set under $f$ is embedded in the part $y$ (notice here how ``inclusion of points in subsets" is changed to embedding of parts). 

This approach has the following problem. It leads to the construction of a set of infinitely many parts of $S$ and, thus, to the undesirable conclusion that the preimage of a shape is a set (of parts), not a shape. A section of this set for the above example is shown graphically below

\begin{figure}[!ht]
\centering
\includegraphics{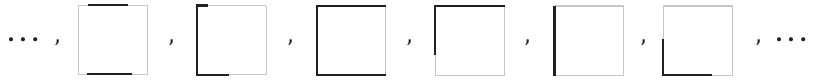}
\label{infinite-set}
\end{figure}

The resolution of this problem, however, is simple enough and it reflects our visual expectations. We will instead define the preimage $f^{-1}(y)$ of a part $y$ of $S'$ as:

\vspace{12pt}
\centerline{ \emph{The largest part of $S$ whose image under $f$ is embedded in $y$}.}
\vspace{12pt}

\noindent In the above example, this largest part is obviously the third part in the set given graphically. Moreover, it is not hard to see that the largest part in this set is equal to the \emph{sum} of all the parts in the set. Although it may be a sum over infinitely many parts, this sum will always result to a shape because the largest part is always a member of this set and covers every other part in it.   

To convert the above definition of the preimage in terms of the ``largest part" into notation, first define the following set

\vspace{12pt}
\centerline{$U_y = \{\; x \;|\; x \leq S \; \textrm{and} \; f(x) \leq y\;\}$}
\vspace{12pt}

\noindent for a part $y$ of the shape $S'$. Then the preimage of part $y$ is the ``maximum" element of $U_y$,

\vspace{12pt}
\centerline{$f^{-1}(y) = \sup U_y$.}
\vspace{12pt}

\noindent This maximum element is the smallest part in $U_y$ greater than any other part in $U_y$. 

Such an element will always exist if the given part $y$ has a preimage. Moreover, if there is at least one nonempty part in the preimage set $U_y$, then there must be infinitely many other parts, including the \emph{empty part} (one can always find infinitely many nonempty parts within a given one). The preimage set $U_y$ may also have one part only, namely the empty one. We thus have that the preimage of a given part $y$ of $S'$ is defined (and it is a shape), if and only if, the set $U_y$ is nonempty; it is not defined, when the set $U_y$ is empty---this, intuitively speaking, mimics an all-or-nothing situation. Consequently, when the set $U_y$ is empty, we will say that the preimage $f^{-1}(y)$ of part $y$ is \emph{undefined}. This is another point where shapes diverge from sets with respect to the preimage operation. For sets, if the preimage set is empty, we say that the preimage is the \emph{empty set} ($\emptyset$). It is not possible to mimic this statement for shapes, however, because the \emph{empty set} is a set, not a shape---the two are of a different object type.

Let us state that for any mapping $f: S \rightarrow S'$ between two shapes, the image operation is required to be \emph{order-preserving}: 

\vspace{12pt}
\centerline{$x \leq x' \implies f(x) \leq f(x')$,}
\vspace{12pt}

\noindent that is to say, it preserves the embedding order of any two parts $x$ and $x'$ in the poset of the shape $S$. By ``f is a mapping between shapes" it is always meant that $f$ is order-preserving (it is possible to define order-reversing mappings between two shapes, but their investigation is out of the scope of this paper).

Let us also state, without proof, that the preimage operation is order-preserving, too. However, the preimage may not be defined for all the parts of a shape $S'$ under a given mapping. There may be parts of $S'$ with undefined preimages. As an example, the preimages of the following three parts are undefined

\begin{figure}[!ht]
\centering
\includegraphics{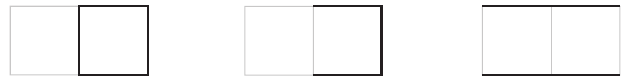}
\label{parts-no-preimage}
\end{figure}

\noindent in the example given previously with the mapping $f(x) = x + A$. And, in fact, there are infinitely many other such parts of $S'$ with undefined preimages.

Now, if the preimage operation is indeed defined for all the parts of $S'$, the following connection between images and preimages can be established for any part $x$ of $S$ and any part $y$ of $S'$:

\vspace{12pt}
\centerline{$f(x) \leq y \iff x \leq f^{-1}(y)$.}
\vspace{12pt}

\begin{proof}
The desired result is immediate from how preimages are defined, and from the assumption that the image and preimage operations are order-preserving. If the image of $x$ is embedded in $y$, then $x$ must be part of $f^{-1}(y)$, which is the largest part of $S$ whose image is embedded in $y$, that is $x \leq f^{-1}(y)$. On the other hand, if $x$ is part of $f^{-1}(y)$ then the image of $x$ must be embedded in $y$, that is $f(x) \leq y$. 
\end{proof}

If the above connection between images and preimages holds for a particular mapping, then (and only then) we obtain the following two convenient embedding relations:

\vspace{12pt}
\centerline{$f(f^{-1}(y)) \leq y$ and $x \leq f^{-1}(f(x))$}
\vspace{12pt}

It is interesting to examine the conditions under which the preimages of all the parts $y$ of $S'$ are defined. One approach is to classify the different types of mappings that may exist between two shapes in terms of \emph{schemas} \cite{Stiny2006}, which are generalized classes of mappings between shapes.



\subsection{Continuous mappings}

A mapping is defined independently of the topologies induced on the two shapes it connects. A mapping only describes how parts of a shape change into parts of another shape. The continuity of the mapping, on the other hand, is a topological issue and it depends on the specific topologies induced on the two shapes. 

In this section, a definition of continuity for shapes is defined that broadly captures the concept of ``structure preservation". A mapping between two shapes is expected to be continuous in a structural sense, when it preserves relations between the open parts in the topologies of the two shapes it connects. Intuitively, we want the open parts of one shape to somehow ``continue" to be related in the same way when a continuous mapping is applied.

The relations between the open parts in a given topology are in principle algebraic and are captured in its lattice of open parts. This highlights the algebraic nature of finite topologies for shapes and provides for a useful definition of continuity.

\begin{figure}[!ht]
\centering
\includegraphics{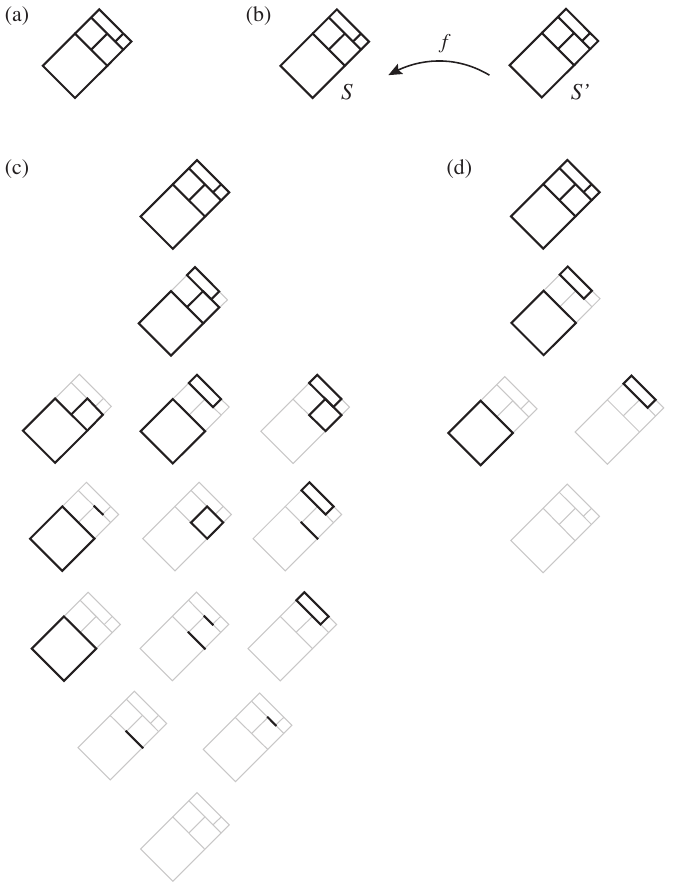}
\caption{Continuity with an identity mapping. (a) A shape $S$, (b) a mapping $f$ from $S$ to itself. $S'$ is the shape $S$ induced with the topology in (d), and $S$ is the same shape induced with the topology in (c).}
\label{Figure21}
\end{figure}

Let $S$ and $S'$ be two shapes and, respectively, $OS$ and $OS'$ be their lattices of open parts. A mapping $f: S \rightarrow S'$ is \emph{continuous} if the mapping $f^*: OS' \rightarrow OS$, defined as $f^*(D) = f^{-1}(D)$, for every open part $D$ in $OS'$, is a \emph{lattice homomorphism}.

Applications of this definition are in \emph{Examples 2a and 2b}. It is implied in the above definition that when $f$ is continuous, the preimage of every open part in $OS'$ is an open part in $OS$ (a result analogous to topological spaces). As a lattice homomorphism, $f^*$ preserves sums and products of any two open parts of $OS'$:

\vspace{12pt}
\centerline{$f^*(C \cdot D) = f^*(C) \cdot f^*(D)$ and $f^*(C + D) = f^*(C) + f^*(D)$.}
\vspace{12pt}

\noindent It follows from (either of) these two equations that the mapping $f^*$ is order-preserving:

\vspace{12pt}
\centerline{If $C \leq D$ in $OS'$ $\implies f^*(C) \leq f^*(D)$ in $OS$.}
\vspace{12pt}

\noindent The proof of this statement is analogous to that for lattices in general (notice that $C \leq D$ implies $C \cdot D = C$ and $C + D = D$).

Two or more parts of the shape $S'$ may have the same preimage under $f$. So the backward mapping $f^*$ is expected to be many-to-one in the general case. Depending on the mapping $f$, however, $f^*$ can be one-to-one, too, in which case we have \emph{structure embedding} (a one-to-one lattice homomorphism). 

By definition, the forward mapping $f$ maps the shape $S$ into the shape $S'$, so that $f(S) = S'$. Thus, $f^*(f(S)) = f^*(S') = S$ (because $S$ is the largest shape embedded in $S'$ under $f$) and it follows that $f^*$ preserves the top element. However, the mapping $f^*$ does not preserve the bottom element, that is, $f^*(0)$ is not necessarily equal to 0 (see \emph{Example 2}). Spatially speaking, the preimage of the \emph{empty shape} may be a nonempty part (What happens with mappings that erase parts?). One example of a mapping that, if continuous, will always preserve both the top and bottom elements is the identity mapping (e.g., \emph{Example 2a}).

The following result is analogous to ``closure under image" in topological spaces for continuous functions. Let $f: S \rightarrow S'$ be a continuous mapping between the shapes $S$ and $S'$. Then,

\vspace{12pt}
\centerline{For every part $x$ of $S$, one has $f(\overline{x}) \leq \overline{f(x)}$.}
\vspace{12pt}

\begin{proof}
Let $D$ = $\overline{f(x)}$ be the closure of $f(x)$ in the topology for $S'$. Since $f$ is continuous by assumption, $f^*(D) = f^{-1}(D)$ is open in the topology for $S$. Then $x$ is a part of $f^{-1}(D)$, because $f^{-1}(D)$ is the largest part of $S$ whose image is embedded in $D$. We thus have $x \leq \overline{x} \leq f^{-1}(D)$ or $f(\overline{x}) \leq f(f^{-1}(D)) \leq D$, so that $f(\overline{x}) \leq \overline{f(x)}$.
\end{proof}

\vspace{12pt}
\noindent \emph{EXAMPLE 2.\;\;\;(a) Let $\mathcal{T}$ and $\mathcal{T}'$ be two comparable topologies for the shape in Figure 21a. Suppose $OS$ and $OS'$ are their lattices of open parts, shown in Figure 21c and 21d, respectively. For ease of reference, call $S$ the shape in Figure 21a induced with topology $\mathcal{T}$ and $S'$ the same shape induced with the topology $\mathcal{T}'$. Define the mapping $f: S' \rightarrow S$ to be the identity map $f(y) = y$, for every part $y$ of $S'$, which is shown diagrammatically in Figure 21b. Mapping $f$ is not continuous. The preimage $f^{-1}(C)$ is not open in $\mathcal{T}'$ for all parts $C$ open in $\mathcal{T}$, and so a homomorphism from $OS$ to $OS'$ is not defined. On the other hand, the inverse mapping of $f$, namely, $f^{-1}: S \rightarrow S'$ is continuous, because $(f^{-1})^{-1}(D) = f(D)$ is open in $\mathcal{T}$, for all parts $D$ open in $\mathcal{T}'$, i.e. a homomorphism is formed from $OS'$ to $OS$. Because the two topologies are defined on the same shape, this implies that every open part in $\mathcal{T}'$ is also an open part in $\mathcal{T}$. Hence, $\mathcal{T}$ is finer than $\mathcal{T}'$.}

\emph{(b) Let $f: S \rightarrow S^{+}$ be the mapping shown in Figure 22a, defined as}

\vspace{12pt}
\centerline{\emph{$f(x) = x - B,$}}
\vspace{12pt}

\noindent \emph{for every part $x$ of $S$, and the symbol $B$ represents the following shape} 

\begin{figure}[!ht]
\centering
\includegraphics{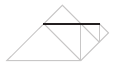}
\label{PartB-sub}
\end{figure}

\noindent \emph{that is subtracted from $S$ to get $S^{+}$. Suppose the shapes $S$ and $S^{+}$ are induced with the topologies in Figure 22b and 22c, respectively. Then, the mapping $f$ is continuous. Moreover, notice that $f^*(S^{+}) = S$ but $f^*(\emph{0}) = B \neq \emph{0}$.}

\begin{figure}[!ht]
\centering
\includegraphics[scale=0.92]{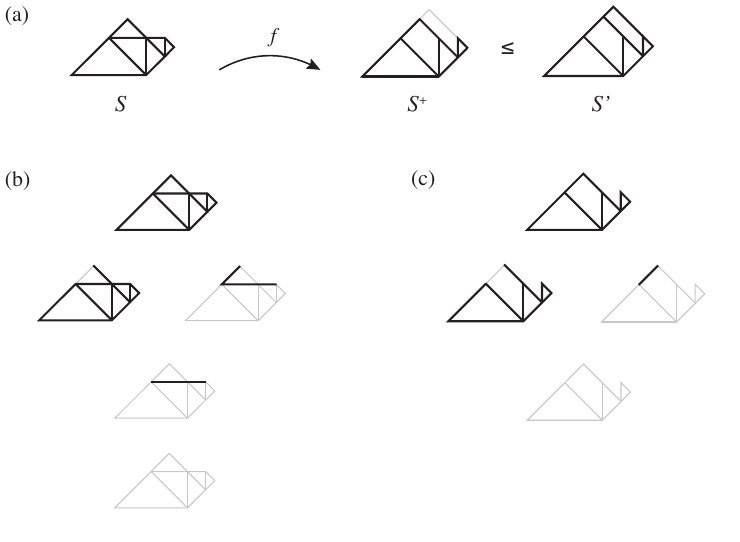}
\caption{(a) An example of a continuous mapping $f$, from a shape $S$ to a proper part $S^+$ of a shape $S'$. The shapes $S$ and $S^+$ possess, respectively, the topologies shown in (b) and (c).}
\label{Figure22}
\end{figure}

\section{Topological connectedness}

\subsection{Visual versus structural connectedness}

The last topic of this paper is topological connectedness. In general topology, connectedness is the study of how topological spaces ``split" into separate pieces, and is perhaps one of the very few topics that have a, more or less, point-free formulation (it is still set-theoretical, but most of the time it can work without mentioning points; see \cite{Munkr00} for the required background). It may seem that shapes can have an exact counterpart, and indeed the basics are as in topological spaces. As we focus, however, on how appearance interacts with connectedness, the results are not as one would expect them to be (at least if one has only point-sets in mind, and not drawings, models and other spatial material).

First, let's define topological connectedness by mimicking the classical version. Let $S$ be a shape and $\mathcal{T}$ a topology for $S$. By \emph{separation} of $S$, define a pair $C$, $D$ of nonempty disjoint open parts in $\mathcal{T}$, such that $C$ + $D$ = $S$. Say that $S$ is \emph{topologically connected}, or just \emph{connected}, if there exists no separation of $S$ in $\mathcal{T}$.

Connectedness, as formulated above, can apply not only to shapes but also to parts. For the subshape topology (\emph{Section 3.3}) associated with a nonempty part of a shape, we have an additional formulation (given here without proof).

Let $S$ be a shape with topology $\mathcal{T}$. If $x$ is a nonempty part of $S$ with subshape topology $\mathcal{T}_x$, a separation of $x$ is a pair of disjoint nonempty parts $A$ and $B$ such that $A$ + $B$ = $x$. The part $x$ is connected, if there exists no separation of $x$ in $\mathcal{T}_x$. 

Note the parts $A$, $B$ that form the separation are open parts in $\mathcal{T}_x$, but don't have to be open in the original topology $\mathcal{T}$ on $S$. The analogous formulation of connectedness for subspaces in general topology, requires the subsets forming the separation to not only be disjoint but to also not share limit points, i.e., each must be disjoint from the closure of the other \cite{Munkr00}. This is automatically granted for disjoint open parts of shapes, because open parts are equal to their closures (\emph{Section 4}).

We can examine some of the topologies given so far. The shape in Figure 14a with the topology in Figure 15a is disconnected. The same shape is connected with the topologies in Figure 16a and 16c, but disconnected with the topology in Figure 16b. The shape in Figure 21a, with either one of the topologies in Figure 21c and 21d, is connected. The shape in the first subshape topology in Figure 15d is connected. The second and third shapes, in the same figure, are disconnected.

Connectedness is a topological (structural) property, since it is determined exclusively from the open parts of a shape relative to a particular topology. Quite like in topological spaces, we have that a shape $S$ is connected, if and only if, it has exactly two complemented open parts, namely, $S$ and \emph{0} (if there is another open part $C$ whose relative complement is open, then $S$ has a separation, namely, $C$ and $S$ - $C$).

Important differences between topological spaces and shapes are revealed when we focus on how appearance interacts with connectedness. 

\begin{figure}[ht]
\centering
\includegraphics[scale=0.83]{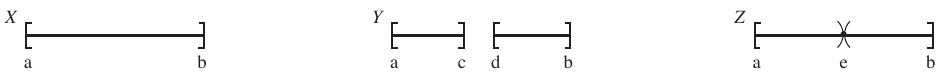}
\caption{Three subspaces of the real line.}
\label{Figure23}
\end{figure}

Consider the three subspaces of the real line $\mathbf{R}$ in Figure 23: $X = [a, b]$, $Y = [a, c] \cup [d, b]$ and $Z = [a, e) \cup (e, b]$. All three subspaces \emph{inherit} an infinite topology from the underlying ``standard" topology for $\mathbf{R}$. $X$ is an interval and, like all intervals of $\mathbf{R}$, it is connected \cite{Munkr00}. $Y$ is the set $X$ with the subset ($c$, $d$) removed. It is disconnected, because it is the union of two disjoint intervals which are open in the subspace topology for $Y$. $Z$ is the set $X$ with the single point $e$ removed. It is also disconnected, because it is the union of two intervals which are open in the subspace topology for $Z$ and neither contains a limit point of the other. Notice how we automatically went from a connected subspace ($X$) to disconnected subspaces ($Y$ and $Z$) by removing pieces. 

With respect to the three subspaces, \emph{structure matches with appearance}: $X$ ``looks connected" (and is structurally connected) because it comes as a single piece; $Y$ ``looks disconnected" (and is structurally disconnected) because it is made of two separate pieces, and the same observation applies to $Z$, too (even though we can't really ``look" at what remains after a single point is removed). Similar scenarios can be devised for other types of shapes, too, for example shapes made of arrangements of planes where each plane is considered as a subspace of $\mathbf{R}^2$.

\begin{figure}[ht]
\centering
\includegraphics[scale=0.83]{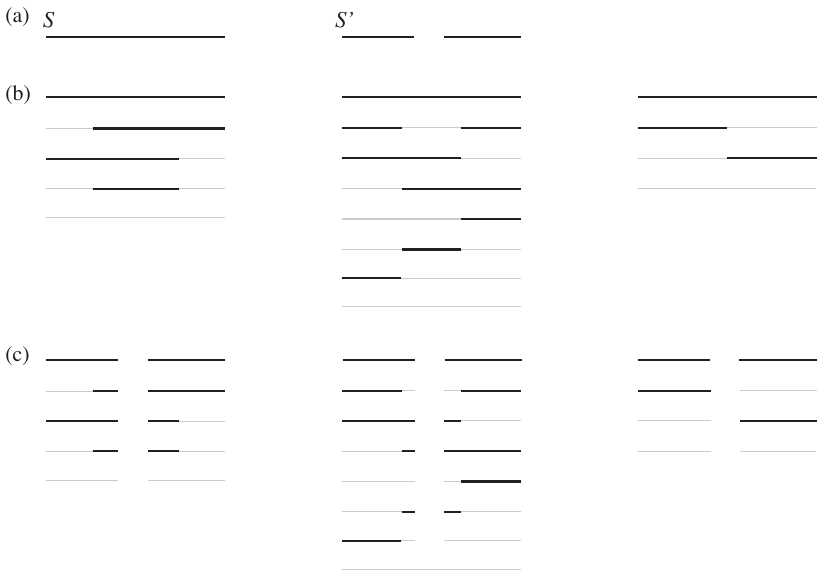}
\caption{(a) A shape $S$ made of a single line and a shape $S'$ obtained by removing a part from $S$. (b) Three topologies for $S$. $S$ is topologically connected only with the first topology. (c) Three subshape topologies for $S'$ (treating $S'$ as a part of $S$). $S'$ is topologically disconnected with the second and third subshape topologies. }
\label{Figure24}
\end{figure}

Unlike (sub)spaces, shapes in algebras $U_i$ do not ``lie" in an ambient space that has a predefined analysis into specific open parts (e.g., subspaces of $\mathbf{R}$ come with a predefined granularity). Shapes in $U_i$ are unanalyzed. Topologies are not inherited but \emph{assigned} through interpretations of their appearance. We can thus have shapes that ``look connected" or that ``look disconnected" induced with either connected or disconnected topological structures, in any number of different ways.

For example, suppose that subspace $X$ is represented as a line in algebra $U_1$, denoted with $S$, as in Figure 24a. As a single line, $S$ ``looks connected". Topologically, however, it can be either connected or disconnected, depending on the open parts recognized in it. Three different finite topologies for $S$ are in Figure 24b; $S$ is connected with the first, but disconnected with the second and third. 

Suppose next that we remove a part from $S$ to obtain the shape $S'$ in Figure 24a (this is subspace $Y$ represented as two lines in algebra $U_1$). That $S'$ now ``looks disconnected"---it comes in two separate pieces---is not an indication (or a consequence) of a disconnected topological structure. $S'$ is a new shape. It can inherit topologies from $S$ through subshape topology constructions, such as the ones in Figure 24c; in this case, $S'$ is still connected with the first subshape topology, and disconnected with the second and third. Alternatively, $S'$ can be induced with completely new topologies, each one with its own connectedness conditions.

Topologies given in previous sections of this paper illustrate it, too. The shape in Figure 14a, for example, ``looks disconnected" because it comes in two pieces. It is topologically connected with respect to the topologies in Figure 16a and 16c, but topologically disconnected with respect to the topologies in Figure 15a and 16b.

The intuitive observation that a shape ``looks connected" or ``looks disconnected" can be expressed in a formally precise way through the concept of \emph{touching} for basic elements in algebras $U_i$, which is defined in \cite[pp. 171-173]{Stiny2006}. Using touching, we can state that a shape \emph{looks connected} (or \emph{is visually connected}) if it is made of basic elements of the same kind that touch\footnote{Touching is defined recursively: basic elements of the same kind \emph{touch} if there are basic elements embedded in each with boundary elements that touch. For example, pairs of lines touch when there are endpoints that do; pairs of planes touch when there are lines embedded in their boundaries with endpoints that do, and so on. One can classify the different ways in which basic elements touch and a demonstration is given in \cite{Stiny2006}.}; otherwise, it \emph{looks disconnected} (or \emph{is visually disconnected}). In Figure 25, shapes are classified according to whether they look connected or disconnected in this visual sense. 

\begin{figure}[ht]
\centering
\includegraphics[scale=0.83]{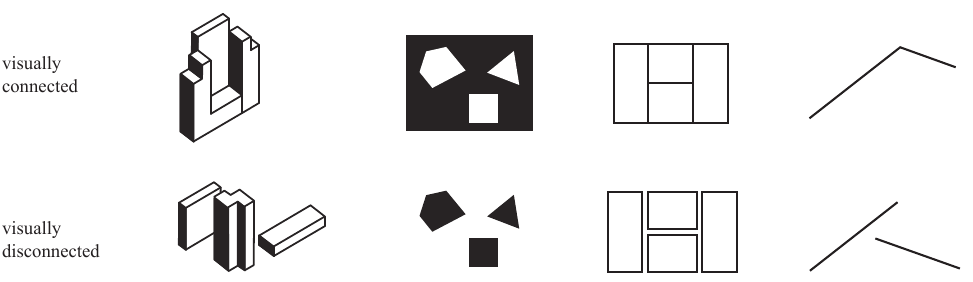}
\caption{Examples of shapes that look connected and shapes that look disconnected.}
\label{Figure25}
\end{figure}

Visual connectedness is about appearance and is independent of structural connectedness. Structural connectedness is a topological property. It depends only on the topology induced on a shape, that is, on the open parts recognized in it (connectedness is not a property to determine just by looking at a shape). As a bottomline:

\vspace{12pt}
\centerline{\emph{Appearance is independent of structure}.}
\vspace{12pt}



\subsection{Further characterizations of connectedness for shapes}

The following result is similar to topological spaces. Let $S$ be a shape with topology $\mathcal{T}$. If a pair $C$, $D$ of open parts forms a separation of shape $S$, and a part $x$ is connected in its subshape topology $\mathcal{T}_x$, then $x$ is part of either $C$ or $D$. 

There is a simple proof analogous to that for spaces, which I omit. Consider an example instead. Suppose $S$ is the line in Figure 24a and is induced with the topology in Figure 26a. The parts $C$ and $D$ in the figure form a separation of $S$. The subshape topology associated with a part $x$ of $S$ falls under the three different cases I show in Figure 26b. If $x$ extends over both $C$ and $D$, its subshape topology will be disconnected. It will be connected, however, if $x$ is entirely embedded in either $C$ or $D$.

\begin{figure}[ht]
\centering
\includegraphics[scale=0.83]{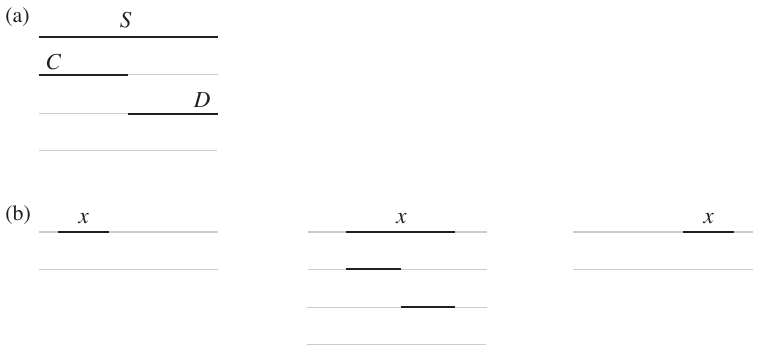}
\caption{(a) A shape $S$, made of a single line, with disconnected topology. Open parts $C$ and $D$ form a separation for $S$. (b) The subshape topology for a part $x$ is connected when $x$ is entirely embedded in either $C$ or $D$.}
\label{Figure26}
\end{figure}

A shape $S$ can be connected with respect to a topology $\mathcal{T}$, but certain open parts of it may be disconnected in their associated subshape topologies. Hence, the following complementary notion.

Say that $S$ is \emph{locally connected at C}, if the open part $C$ is connected in its associated subshape topology $\mathcal{T}_C$. If $S$ is locally connected at every open part $C$, then $S$ is said to be \emph{locally connected}. 

If a shape $S$ with topology $\mathcal{T}$ is locally connected, then it is connected (a trivial observation); the opposite, however, is obviously not true. For example, the shape in Figure 21a, with either one of the topologies in Figure 21c and 21d, is connected but is not locally connected at certain open parts. 

It is reasonable to ask whether we can have a finite topology for a shape in which \emph{every nonempty open part} is disconnected in its associated subshape topology---a ``totally disconnected" finite topology, so to speak. Since topologies for shapes are finite, and there are no points, open parts cannot be partitioned indefinitely in them. Ultimately, there must be certain ``minimal" open parts which are connected in their subshape topologies, even if every other open part is disconnected. These minimal open parts are precisely the elements in a reduced basis. The following finitistic definition of a totally disconnected topology for a shape derives from this fact.

A shape $S$ with topology $\mathcal{T}$ is \emph{totally disconnected} if $S$ is disconnected and the only connected open parts in $\mathcal{T}$ are the basis elements in the reduced basis of $\mathcal{T}$.

An example of a topology that makes the shape in Figure 14a totally disconnected (in this finitistic sense), is shown in Figure 27. The reduced basis of this topology is the set of four nonempty open parts in the bottom row.

\begin{figure}[ht]
\centering
\includegraphics{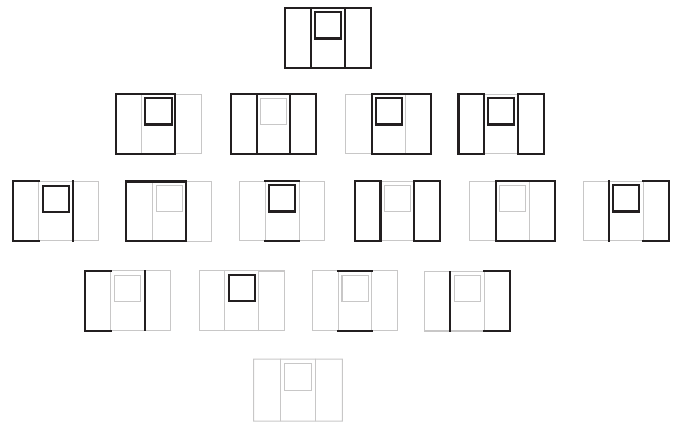}
\caption{A totally disconnected topology for the shape in Figure 14a.}
\label{Figure27}
\end{figure}

Let the shape $S$ with topology $\mathcal{T}$ be totally disconnected, and let $\mathcal{B}$ be the reduced basis of $\mathcal{T}$. The following are equivalent statements:

\begin{enumerate}
\item[(1)] $S$ is totally disconnected.

\item[(2)] $\mathcal{B}$ consists of disjoint nonempty parts only (excluding the \emph{empty shape}).

\item[(3)] $\mathcal{T}$ consists of closed-open parts only.

\item[(4)] $\mathcal{T}$ determines a finite Boolean algebra over the closed-open parts of $S$, with bottom element \emph{0} and top element $S$ itself.
\end{enumerate}

\begin{proof}
The proofs of implications (1) through (4) are almost automatic.

(1)$\implies$(2): Suppose for contradiction that $b_i$ and $b_j$ are two basis elements in $\mathcal{B}$ that are not disjoint. Then, there must be an open part $C$ formed as the sum of $b_i$ and $b_j$, which is connected. But we already know from assumption (1) that such an open part cannot exist. 

(2)$\implies$(3): By definition, $S$ can be expressed as a sum of disjoint nonempty basis elements $b_1, ..., b_n$ of $\mathcal{B}$. Given $C$, an open part in $\mathcal{T}$, $C$ is either a basis element or a sum of basis elements, i.e. $C = b_1 + ... + b_k$ where $1 \leq k \leq n$. Then, rewrite $S = C + b_{k+1} + ... + b_n$, so that $S - C = b_{k+1} + ... + b_n$. The right-hand side of the latter equality is an open part in $\mathcal{T}$, by definition. Thus, we have shown that $C$ and its relative complement $S - C$ are open in $\mathcal{T}$, and therefore $C$ is closed-open.

(3)$\implies$(4): $\mathcal{T}$ is, by definition, a finite, distributive lattice, with bottom element \emph{0} and top element $S$. By assumption (3), every member of $\mathcal{T}$ is closed-open and thus comes with a relative complement. Hence, $\mathcal{T}$ is a Boolean algebra over the closed-open parts of $S$.

(4)$\implies$(1): From Stone duality, the finite Boolean algebra determined from assumption (4) corresponds to a totally disconnected space (a Stone space). The isolated points of this space are in bijection with the atoms of the Boolean algebra, and the closed-open subsets are formed by unions of those isolated points. As shown in \emph{Section 5}, this space is isomorphic to a finite shape topology.
\end{proof}

Finally, statement (2) given above implies that we can form a totally disconnected topology for a shape from a basis consisting of disjoint nonempty parts. In how many ways can this happen? Indefinitely many. Any shape can be decomposed into finitely many disjoint parts (that sum to the shape), in any number of different ways, each one resulting in a different totally disconnected topology (contrast this result with finite spaces, where totally disconnected topologies are formed uniquely).

\section{Discussion}

\subsection{More topological concepts for shapes}

\emph{Sections 3} through \emph{7} are not exhaustive, but they cover a wide range of foundational topological concepts. Many other concepts are possible though.

One important topic in topology that is left without mention are the so-called \emph{countability} and \emph{separation axioms}. Traditionally, these axioms reflect properties of the overall topological structure assigned to an object and are more often used in the study of infinite objects rather than finite ones (consult \cite{Kurat72} and \cite{Munkr00} for the required background). It is possible to mimic some of the countability axioms in the context of shapes using the constructions presented in this paper (e.g., \emph{Section 3}). For example, any shape with any finite topology is \emph{second-countable} because the topology is generated by a countable basis.

The separation axioms characterize the ways that the members of a topology are ``separated" from one another. In general, they are defined in terms of points and so a direct transfer to shapes is not always possible or intuitive. For example, the T$_0$ separation axiom is not really applicable to shapes in algebras $U_i$, for $i > 0$, because it refers to the separation of points. It is applicable, however, to shapes made with points in algebra $U_0$, as discussed in \cite{Haridis2020}. To formulate separation axioms, one would need to cast them in terms of embedding relations between open parts.

Other topics I have not discussed but that are plausible to pursue for finite topologies on shapes, are the constructions of \emph{product topologies}, \emph{homeomorphisms}, \emph{neighborhood systems}, and \emph{metrizability}. The definition of a continuous mapping in \emph{Section 6.2} can be used to obtain a notion of homeomorphism between two shapes. Then one could ask questions of the sort: Are a square and a triangle, as shapes in algebra $U_1$, homeomorphic to each other? How about a square and a circle? With respect to metrizability, under what conditions are finite topologies for shapes in algebras $U_i$ metrizable? Is there something meaningful to say about the relationship between metrizable topologies for shapes and the classical metric spaces in general topology? In light of a point-free framework for topology, one could also ask what happens with the classical concept of convexity, which is normally understood in point-set theoretical terms.

It would be also interesting to see what extensions are required to the material developed in this paper in order to accommodate topology on shapes equipped with weights (e.g., colors, thicknesses, and so on) or shapes made with basic elements from composite algebras. For example, how would a topology be defined when a shape is made up of both points and lines or both lines and planes? Last, it is worth exploring infinite topologies for shapes in algebras $U_i$, for $i > 0$, to see what topological constructions of this paper transfer smoothly in the infinite case. 

\subsection{Mathematics of shapes (in art and design)}

The abstract idea of ``structure" is a recurring theme in the mathematical description and analysis of design objects in architecture, visual arts, city planning and other design areas. There are many ways to talk about structure and perhaps the only assumption common to all is a satisfactory answer to the question ``What is it [the object] made of?" \cite{Runci69}. When asking about the structure of an object, it is important to have the following distinction in mind, which recapitulates the motivation for the work I presented in this paper.

If the object under consideration is analyzed, that is, when the elements it is composed of are given separately and individually (e.g., as the points of a topological space), to determine its structure, one is led to restrict attention to the interconnections and interdependencies of the individual elements (for example, topology with finite sets of points, or with the vertices of a graph, falls under this category). What happens in the case where the object under consideration comes without a given, or apparent, analysis into separate and individual elements? This is exactly the case of shapes in algebras $U_i$ for $i > 0$ (and, by analogy, the case of drawings, pictures, models, in art and design). Such shapes are inherently \emph{unanalyzed} in the sense that they do not come equipped with a definite subdivision into parts. For any given shape, to say what the elements are ``upon which structure can be induced and thereafter studied mathematically", an interpretative act must be introduced. This interpretative act may take place not singly but multiply: shapes (made of either lines, planes or solids) have finite definition, but are indefinitely interpretable into parts (\emph{Section 2}).

This distinction is important for another reason. Most work on structural mathematics in design, especially in the fields of architecture and planning, is about the study of relations on fixed and definite parts of a given object---the ``analyzed approach" I described just above. There are multiple examples from the literature; see \cite{Alexan64, AlexanPoyn67R3, Atkin74} and \cite{MarchStead74, Steadman83}. While it is a reasonable approach, it is perhaps only because no domain of mathematics that deals explicitly with structure (e.g., finite combinatorial topology, graph theory, combinatorics) has anything to do with appearance. It is thus hoped that this paper brings into awareness that the structure (topology in this case) of objects in art and design can indeed be studied mathematically, but yet in a way that retains their pictorial/spatial characteristics and at the same time takes into consideration their aesthetic, interpretative capacity.

By ``mathematics of shapes" then (algebra or topology, for instance), we do not only mean the application of mathematics by injecting shapes with number or point-set theoretical properties and patterns. More importantly, we mean the mathematics that arises naturally when we take the pictorial/spatial characteristics of shapes as the starting point of a mathematical investigation.

\section*{Acknowledgement(s)}

This paper is a preprint version of a journal article published in the \emph{Journal of Mathematics and the Arts} (Taylor $\&$ Francis), available at: \url{https://doi.org/10.1080/17513472.2020.1723828}{}.

\noindent \textbf{Cite as}:

\noindent Haridis A. (2020) Structure from appearance: topology with shapes, without points. \emph{Journal of Mathematics and the Arts 14}(3): 199-238.


\section*{Declaration of interest statement}

No potential conflict of interest was reported by the authors.

\bigskip

\section{Appendices}

\appendix
\section{Proofs related to the basis sets $\mathcal{B}$ and $\mathcal{B}_x$.}

\subsection{} \noindent Let $S$ be a shape and $\mathcal{B}$ a nonempty set of parts of $S$, satisfying conditions (1)-(3) of \emph{Section 3.2}. Then $\mathcal{B}$ is the \emph{unique minimal basis} (or the \emph{reduced basis}) for a topology $\mathcal{T}$ of $S$.

\begin{proof} Conditions (1) and (2) are basic conditions for a basis that $\mathcal{B}$ satisfies. Choose $b$, any element of $\mathcal{B}$ and suppose, for contradiction, that the set $\mathcal{B} \setminus \{b\}$ still describes $b$. That is, there are two other elements in $\mathcal{B} \setminus \{b\}$ that sum to $b$. But we know that such elements cannot exist because, by condition (3), $b$ must be one of the elements of which it is the sum. It follows that $b$ must be contained in $\mathcal{B}$ and thus in any other basis $\mathcal{B}'$ for the same topology. Since $b$ is arbitrary, we conclude that $\mathcal{B} \subset \mathcal{B}'$.

That $\mathcal{B}$ is unique follows immediately by the preceding argument. If $\mathcal{B}'$ were in fact a second minimal basis, then $\mathcal{B}' \subset \mathcal{B}$, so that $\mathcal{B} = \mathcal{B}'$, as desired. \end{proof}

\subsection{} 

\noindent Let $\mathcal{B}$ be a basis for a topology $\mathcal{T}$ on $S$ and $x$ a part of $S$. Then, the set $\mathcal{B}_x$, defined in \emph{Section 3.4}, is a basis for the subshape topology $\mathcal{T}_x$.

\begin{proof} We first show that the set $\mathcal{B}_x$ is a basis for a topology on $x$, by showing that it satisfies conditions (1) and (2) given in \emph{Section 3.2}.

The set $\mathcal{B}_x$ satisfies condition (1) because it is formed by products between $x$ and every basis elements in $\mathcal{B}$: 

\vspace{0.1in}
\centerline{$b_1 \cdot x + ... + b_n \cdot x$ = $(b_1 + ... + b_n)\cdot x = S \cdot x = x$}
\vspace{0.1in}

\noindent where $b_1, ..., b_n$ are the basis elements of $\mathcal{B}$. 

To show that $\mathcal{B}_x$ satisfies condition (2), take any two basis elements $b_1 \cdot x$ and $b_2 \cdot x$ of $\mathcal{B}_x$. We have 

\vspace{0.1in}\centerline{ $(b_1 \cdot x) \cdot (b_2 \cdot x) = (b_1 \cdot b_2) \cdot x $}
\vspace{0.1in}

\noindent where the product $b_1 \cdot b_2$ is (by definition) either equal to a basis element $b_3$, or a sum $\sum b_i$ of basis elements $b_i$ of $\mathcal{B}$. If it is the former, $b_3 \cdot x$ is in $\mathcal{B}_x$, by definition. If it is the latter, $(\sum b_i) \cdot x = \sum (b_i \cdot x) $ is a sum of basis elements of $\mathcal{B}_x$, as desired.

We have shown that the set $\mathcal{B}_x$ is a basis for a topology on $x$; call this topology $\mathcal{T}'$. We now show that the topology $\mathcal{T}'$ is indeed equal to the subshape topology $\mathcal{T}_x$. 

Take $C \cdot x$ any open part in $\mathcal{T}_x$. Then, $C \cdot x$ can be written as 

\vspace{0.1in}\centerline{$(b_1 + ... + b_k) \cdot x = b_1 \cdot x + ... + b_k \cdot x$}
\vspace{0.1in}

\noindent for basis elements $b_1, ..., b_k$ of $\mathcal{B}$ that sum to $C$. Each member $b_i \cdot x$, for all $i = 1, ..., k$, of the product in the right hand side is a basis element in $\mathcal{B}_x$, by definition. Therefore, $C \cdot x$ is a sum of basis elements of $\mathcal{B}_x$, which we know is open in $\mathcal{T}'$. Conversely, if $D$ is an open part in $\mathcal{T}'$, then $D$ is the sum of basis elements of $\mathcal{B}_x$. Since each of those basis elements belongs to $\mathcal{T}_x$ and $\mathcal{T}_x$ is closed under sum, $D$ belongs in $\mathcal{T}_x$.
\end{proof}



\appendix

\end{document}